\documentclass[12pt,a4paper]{scrartcl}
\usepackage{amsmath}\usepackage{amssymb}
\usepackage{amsthm}
\usepackage[english]{babel}\usepackage{german}
\usepackage{nicefrac} \usepackage{floatflt}
\usepackage[pdftex]{graphicx}  
\usepackage[maxbibnames=10]{biblatex}
\newcommand{\R}{\mathbb{R}}
\newcommand{\Q}{\mathbb{Q}}
\newcommand{\Z}{\mathbb{Z}}

\newcommand{\tr}{\mathrm{tr}}

\newcommand{\C}{\mathbb{C}}
\newcommand{\CC}{\mathcal{C}}

\newcommand{\OK}{\mathcal{O}_K}

\newcommand{\OKnu}{\mathcal{O}_{K_\nu}}

\newcommand{\Al}{A_\lambda}
\newcommand{\Sf}{S_<^{\mathrm{fin}}}
\newcommand{\Si}{S_<^{\infty}}
\newcommand{\Alf}{\Al^{\mathrm{fin}}}
\newcommand{\Ali}{A_\lambda^\infty}

\newcommand{\Ab}{A_\beta}

\newcommand{\EE}{\mathbb{E}}
\newcommand{\PP}{\mathbb{P}}

\newcommand{\N}{\mathbb{N}}
\newcommand{\mo}{\; \mathrm{mod} \;}
\newcommand{\eps}{\varepsilon}

\makeatletter
\newtheorem*{rep@theorem}{\rep@title}
\newcommand{\newreptheorem}[2]{%
\newenvironment{rep#1}[1]{%
 \def\rep@title{#2 \ref{##1}}%
 \begin{rep@theorem}}%
 {\end{rep@theorem}}}
\makeatother

\title{On absolute continuity and maximal Garsia entropy for self-similar measures with algebraic contraction ratio }
\author{ Lauritz Streck \thanks{Institution: University of Cambridge.  Address: Centre for Mathematical Sciences, Wilberforce Road, CB3 0WA Cambridge.  E-mail: ls909@cam.ac.uk}}
\date{\today}

\begin{document}

\newtheorem{theorem}{Theorem}[section]
\newreptheorem{theorem}{Theorem}
\newtheorem{definition}[theorem]{Definition}
\newtheorem{remark}[theorem]{Remark}
\newtheorem{lemma}[theorem]{Lemma}
\newtheorem*{lemma*}{Lemma}
\newreptheorem{lemma}{Lemma}
\newtheorem{corollary}[theorem]{Corollary}
\newtheorem{prop}[theorem]{Proposition}
\newreptheorem{prop}{Proposition}
\newtheorem{claim}[theorem]{Claim}
\newreptheorem{claim}{Claim}
\newtheorem{conj}[theorem]{Conjecture}

\newtheorem{example}[theorem]{Example}

\maketitle

\begin{abstract}
We consider the self-similar measure $\nu_\lambda=\mathrm{law}\left(\sum_{j \geq 0} \xi_j \lambda^j\right)$ on $\R$, where $|\lambda|<1$ and the $\xi_j \sim \nu$ are independent, identically distributed with respect to a measure $\nu$ finitely supported on $\Z$.  One example of this are Bernoulli convolutions. It is known that for certain combinations of algebraic $\lambda$ and $\nu$ uniform on an interval, $\nu_\lambda$ is absolutely continuous and its Fourier transform has power decay; in the proof, it is exploited that for these combinations, a quantity called the Garsia entropy $h_{\lambda}(\nu)$ is maximal. 

In this paper, we show that the phenomenon of $h_{\lambda}(\nu)$ being maximal is equivalent to absolute continuity of a self-affine measure $\mu_\lambda$, which is naturally associated to $\lambda$ and projects onto $\nu_\lambda$.  We also classify all combinations for which this phenomenon occurs: We find that if an algebraic $\lambda$ without a Galois conjugate of modulus exactly one has a $\nu$ such that $h_{\lambda}(\nu)$ is maximal, then all Galois conjugates of $\lambda$ must be smaller in modulus than one and $\nu$ must satisfy a certain finite set of linear equations in terms of $\lambda$.  Lastly, we show that in this case, the measure $\mu_\lambda$ is not only absolutely continuous but also has power Fourier decay, which implies the same for $\nu_\lambda$.   
\end{abstract}

\tableofcontents
\section{Introduction}
\subsection{Motivation} \label{subsec::motivation}
\paragraph{Problem Statement}
Let $\lambda$ be a real number with $|\lambda|<1$ and let $\nu$ be a measure with finite support on the integers. Let $\xi_0, \xi_1, \dots$ be independent, identically distributed with respect to $\nu$.  Define the self similar measure $\nu_\lambda$ on $\R$ by
\[
\nu_\lambda:=\mathrm{law}\left(\sum_{j=0}^\infty \xi_j \lambda^j \right).
\]
Questions about the Hausdorff dimension, absolute continuity and Fourier decay of self-similar measures have a long and varied history.  
We give an overview of some of the results after some definitions. 

For the study of algebraic $\lambda$, a crucial quantity is the \textit{Garsia Entropy} 
\[
h_\lambda(\nu):=\lim_{n \to \infty} \frac{H(X_n)}{n},
\]
where 
\[
X_n:=\sum_{j=0}^{n-1} \xi_j \lambda^j
\]
and $H(X)=-\sum_x \PP(X=x) \log \PP(X=x)$ is the Shannon entropy of a finite random variable.  The limit exists because the sequence $H(X_n)$ is sub-additive, which implies that $\frac{H(X_n)}{n}$ is monotonously decreasing.

An aspect of the complexity of an algebraic $\lambda$ is quantified by its \textit{Mahler measure} 
\[
M_\lambda:=|a_d| \prod_{\sigma: \Q(\lambda) \to \C, |\sigma(\lambda)|>1} |\sigma(\lambda)|,
\]
where $a_d$ is the leading coefficient of the minimal polynomial of $\lambda$ over $\Z[X]$ and $\sigma(\lambda)$ ranges over the Galois conjugates of $\lambda$. 
The Mahler measure is related very naturally to the Garsia entropy because, as is well-known,
\begin{equation} \label{eq::entroineq}
h_\lambda(\nu) \leq \log M_\lambda
\end{equation}
for any $\nu$ (a proof is given in Lemma \ref{lem::separation}).  It is also known that, roughly speaking, $h_{\lambda}(\nu) \to \log M_\lambda$ as $\nu$ expands (for example for $\lambda$ an algebraic unit and $\nu=\mathrm{Unif}\{0, \dots, N\}$ with $N \to \infty$,  \cite{peters}). We say that a pair $\lambda, \nu$ has \textit{maximal entropy} if $h_\lambda(\nu)=\log M_\lambda$.  Classifying all combinations of $\lambda$ and $\nu$ with maximal entropy and deducing properties like absolute continuity and power Fourier decay for such combinations is the aim of this paper.

\paragraph{The role of Garsia entropy in fractal geometry}
Before we say more about maximal entropy,  we explain why the Garsia entropy plays such an important role in the field of fractal geometry.  In 2013,  Hochman published a breakthrough paper in which he used entropy rates to show properties of the measures $\nu_\lambda$ \cite{hochman}.  As a corollary of this, one can deduce that for algebraic $\lambda$,
\[
\mathrm{dim}(\nu_\lambda)=\min\left(1, \frac{h_\lambda(\nu)}{\log \lambda^{-1}} \right),
\]
where $\mathrm{dim}$ denotes the Hausdorff dimension. Using Hochman's result,  Varj\'u showed  that $\mathrm{dim}(\nu^{\mathrm{Ber}}_\lambda)=1$ for transcendental $\lambda \in (0.5,1)$ \cite{varju}.  The best resources to study the respective proofs are \cite{hochmanicm} and \cite{varjuicm}.

Varj\'u's result resolved the \textit{exact overlaps conjecture} on $\R$ for the case of two maps,  corresponding to the two points $\nu^{\mathrm{Ber}}$ is supported on. The exact overlaps conjecture is commonly considered to be one of the most fundamental open problems concerning self-similar measures.  To keep this introduction from getting too technical, we do not state it here but mention that in the statement of the full exact overlaps conjecture,  the measure $\nu^{\mathrm{Ber}}$ is replaced by a probability measure supported on $n$ points.  

Varj\'u's proof uses a result of Breuillard and Varj\'u to the effect that transcendental $\lambda$ with $\mathrm{dim}(\nu_\lambda)<1$ can be approximated arbitrarily well by algebraic $\lambda_n$ with $h_{\lambda_n}(\nu)<\log \lambda^{-1} -\eps$, with $\eps$ depending on $\lambda$ \cite{bv2}.  A crucial part of the proof is to bound the Mahler measure of those $\lambda_n$, which is possible due to the result of Breulliard and Varj\'u \cite{bv} that  for every $\eps$ there is $C=C(\eps)$ such that
\begin{equation} \label{bvresult}
 h_\lambda(\nu^{\mathrm{Ber}}) \leq \log 2-\eps \quad \Longrightarrow \quad M_\lambda \leq C.
\end{equation}

Rapaport and Varj\'u \cite{rapavarju} worked on generalising the result in \cite{varju} from $\nu^{\mathrm{Ber}}$ to a measure supported on three points and made progress, but could not solve this problem completely.  The major roadblock preventing this is the missing generalisation of (\ref{bvresult}) to this setting.  While there are more problems to be resolved for the case of $n$ instead of three points,  Varj\'u's proof does provide a road map, of which a crucial step is once again a suitable generalisation of (\ref{bvresult}).  Improving the understanding of the relationship between the Garsia entropy and the Mahler measure enough could thus potentially lead to major progress on the exact overlaps conjecture.

\paragraph{Computation of the Garsia entropy}
We outline some known methods to estimate the Garsia entropy.  As we will point out at the end of this Subsection, part of our work is to extend the methods sketched here from the Archimedean to the non-Archimedean completions, which we make precise later.

As we recall, the Garsia entropy is defined in terms of the sub-additive sequence $H(X_n)$, with $\frac{1}{n}H(X_n)$ converging to $h_\lambda(\nu)$ from above.  To find lower bounds on the entropy, an important idea is to find a super-additive sequence $L_n$ such that $\frac{1}{n} L_n$ converges to $h_\lambda(\nu)$ from below.  One can then deduce lower bounds on the Garsia entropy by estimating the terms of this sequence.

This idea is a crucial step in the proof of (\ref{bvresult}); the proof operates by finding an appropriate super-additive sequence and bounding its first term from below. This also allows to calculate $h_\lambda(\nu)$ to arbitrary precision (in principle),  provided that $\lambda$ is an algebraic unit. Here,  an algebraic unit is a number such that both $\lambda$ and $\lambda^{-1}$ are algebraic integers. A different super-additive sequence is used in \cite{4author} to give an explicit algorithm to compute $h_\lambda(\nu)$ to arbitrary precision, again provided that $\lambda$ is an algebraic unit. Neither of these methods can be used if $\lambda$ is not an algebraic unit,  like for $\lambda=\frac{2}{3}$.  

In both papers,  one considers a more general self-affine measure,  which lives on some $\R^s$ and is related to the Galois conjugates $\sigma(\lambda)$ of $\lambda$.  To find the super-additive sequence,  one then smooths this measure with an object of proportion $|\sigma(\lambda)|^n$ in each coordinate.   


\paragraph{Known results about absolute continuity and Fourier decay} In this paragraph,  we sketch the examples of $\lambda$ for which the behaviour of $\nu_\lambda$ is known. 
Many of these examples are related to the phenomenon of maximal entropy in some way.   

The origins of the study of Bernoulli convolutions go back to Erd\H os,  who showed in 1939 that for certain algebraic $\lambda \in (0.5, 1)$, $\nu^\mathrm{Ber}_\lambda$ is not absolutely continuous.  Those algebraic numbers are $\lambda$ such that $\lambda^{-1}$ is a Pisot number, which means that $|\sigma(\lambda)|>1$ for all Galois conjugates of $\lambda$.  In \cite{garsia1}, Garsia strengthened Erd\H os' result to show that $\mathrm{dim}(\nu^\mathrm{Ber}_\lambda)<1$ in this case.  To this end, he introduced the Garsia entropy.  He started by showing that in the Pisot case,  $\mathrm{dim}(\nu^\mathrm{Ber}_\lambda)=1$ would imply that the Garsia entropy is maximal.  He then showed that maximal Garsia entropy would imply absolute continuity, in contradiction to the result of Erd\H os. To this day,  Pisot numbers are the only $\lambda \in (0.5,1)$ known for which $\nu^{\mathrm{Ber}}_\lambda$ is not absolutely continuous.   

We now sketch known generic results. They are unrelated to maximal entropy but are given for context. In the range of $\lambda \in (0, 0.5)$,  $\nu^{\mathrm{Ber}}_\lambda$ is a Cantor measure, not absolutely continuous but well understood.  It was eventually shown in the nineties that $\nu^{\mathrm{Ber}}_\lambda$ is absolutely continuous except for a Lebesgue null set of $\lambda$ \cite{solomyak}. This was later strengthened to an exceptional set of Hausdorff dimension zero \cite{shmerkin}, using Hochman's breakthrough cited above.  It is also known that almost every $\lambda$ has power Fourier decay \cite{erdos2},  \cite{60years}. Here, \textit{power Fourier decay} means that there are $\delta, C$ such that 
\[
|\widehat{\nu_\lambda}(v)| \leq C |v|^{-\delta}
\]
for all $v \in \R$, where $\widehat{\nu_\lambda}(v):=\int \mathrm{e}^{2 \pi i vy} d\nu_\lambda(y)$ is the Fourier transform.  

However, for all the strong positive results for almost all $\lambda$, there are very few concrete examples of $\lambda$ for which absolute continuity and power Fourier decay are known. 
One family of concrete examples was given by Garsia \cite{garsia2}, who showed that for $\lambda$ the inverse of a Garsia number $\beta=\lambda^{-1}$, the measure $\nu^\mathrm{Ber}_\lambda$ is absolutely continuous with bounded density.  Here, a Garsia number $\beta$ is an algebraic $\beta$ such that $|\sigma(\beta)|>1$ for all Galois conjugates of $\beta$ and $M_\beta=2$.   Inverses of Garsia numbers are the only $\lambda$ for which $\nu^{\mathrm{Ber}}_\lambda$ attains maximal entropy.  In fact, the natural generalisation of Garsia's ideas exactly is that  maximal entropy implies absolute continuity for all $\lambda$ (provided that $\lambda$ has no Galois conjugates on the unit circle,) as we show in Section \ref{sec::computation}.

It was later shown in \cite{feng} that this example can be generalised to any $\lambda$ with $|\sigma(\lambda)|<1$ for all $\sigma$ and $\nu=\mathrm{Unif}\{0, \dots, M_\lambda-1\}$.  The assumption on $\lambda$ implies that $M=M_\lambda$ is an integer.  Moreover, it was shown that for these examples, $\nu_\lambda$ has power Fourier decay.  The proof of power Fourier decay uses Garsia's ideas to establish absolute continuity of $\nu_\lambda$ as a starting point.  Once again,  for those combinations of $\lambda$ and $\nu$ maximal entropy is attained. 

Even though these are somewhat special examples and we know that absolute continuity and power decay hold generically, there are few other concrete examples.  In the case of absolute continuity,  Varj\'u \cite{bv} showed that there is a very small $c$ such that for all algebraic $\lambda$ with $\lambda>1-c \min\left\{\log M_\lambda, (\log M_\lambda)^{-1.1} \right\}$, the measure $\nu^{\mathrm{Ber}}_\lambda$ is absolutely continuous. Kittle \cite{kittle}  generalised the ideas to show that for algebraic $\lambda$ with Mahler measure close to two, the measure is also absolutely continuous.  In the case of power Fourier decay, our knowledge is even worse: The examples in \cite{garsia1} and \cite{feng} are the only known examples.

\paragraph{Contribution of this paper}
Given that the phenomenon of maximal entropy provides a powerful method to find $\lambda$ and $\nu$ such that $\nu_\lambda$ is absolutely continuous and has power Fourier decay and many of the known examples employ it,  one can ask how special this phenomenon really is.  We answer this question in this paper by classifying the pairs of algebraic $\lambda$ (without conjugate on the unit circle) and measures $\nu$ for which maximal entropy is attained.  

One can also ask about the nature of the connection between absolute continuity and maximal entropy. Here, we show that maximal entropy is in fact equivalent to absolute continuity of a self-affine measure $\mu_\lambda$ projecting onto $\nu_\lambda$ (see Corollary \ref{cor::equiv}).  This paper thus unifies all of the papers (\cite{feng}, \cite{erdos}, \cite{garsia2}, \cite{garsia1}) mentioned in the preceding part.

We briefly sketch the setting of the measure $\mu_\lambda$.  It lives on a space we call  the \textit{contraction space} and denote by $\Al$.  This space is naturally associated to $\lambda$ and is a product of $\R^s$ with some non-Archimedean fields.  The measure $\mu_\lambda$  for $\lambda$ an algebraic unit is the one considered in the papers \cite{bv}, \cite{4author} to calculate the Garsia entropy (see ``Computation of the Garsia entropy''); their method works exactly when $\Al \cong \R^s$ and no non-Archimedean fields are present.  We introduce both $\mu_\lambda$ and $\Al$ formally in Subsection \ref{subsec::concept}.  

In Section \ref{sec::computation},  we develop the non-Archimedean analogue of the method in \cite{bv}, \cite{4author} to estimate the entropy in the non-Archimedean completions.  We prove in Proposition \ref{prop::superadd} that if $\Al$ consists only of non-Archimedean fields,  then there is an analogous super-additive sequence converging to the entropy.  In this setting,  the super-additive sequence is even of a particularly easy form: One has
\[
L_n=H(Y_n),
\]
where $H(Y_n)$ is the Shannon entropy of a random variable taking finitely many values.  Combining this with one of the aforementioned methods for the Archimedean completions,  one can for example calculate $h_\lambda(\nu)$ for any combination of $\nu$  and $\lambda$ (without Galois conjugate on the unit circle)  to arbitrary precision.  

Calculating $h_\lambda(\nu)$ for $\lambda$ which is not an algebraic unit has not been studied much before. This is chiefly because the main interest traditionally has been to study the Bernoulli measure $\nu^{\mathrm{Ber}}$. For a measure supported on only two points, our method is not needed to estimate the Garsia entropy because $X_n$ generates a free random walk unless $\lambda$ is an algebraic unit.  It is only when one considers measures supported on more than two points, as one has to in order to progress on the exact overlaps conjecture,  that our method becomes relevant.

It will become noticeable throughout this paper how much easier proofs often become when one considers the space $\Al$ instead of $\R$. An example of this is the proof of power decay in Theorem \ref{main2}. Because of this, it seems likely that the legwork we do in Sections \ref{sec::abscont}, \ref{sec::computation} and \ref{sec::notalgint} to systematically explore and prove properties of $\Al$ and $\mu_\lambda$ will be valuable to future mathematicians working on understanding the Garsia entropy.

\subsection{Results}
For an algebraic number $\lambda$ of degree $d$, there are $d$ field embeddings $\sigma: \Q(\lambda) \to \C$. The algebraic numbers $\sigma(\lambda)$ for $\sigma$ a field embedding to $\C$ are called the \textit{Galois conjugates} of $\lambda$. 
Our first result is a negative result when the phenomenon of maximal entropy can not happen. 
\begin{theorem} \label{main1}
Let $\lambda \in \C$, $|\lambda|<1$, be an algebraic number without a Galois conjugate on the unit circle. Let $\nu$ be a finitely supported probability measure on $\Z$. If there is a field embedding $\sigma: \Q(\lambda) \to \C$ such that $|\sigma(\lambda)|>1$, then $h_\lambda(\nu)<\log M_\lambda$.
\end{theorem}
Br\'emont \cite{bremont} considered a similar question in a more general case (when the contraction ratios can be powers of the same $\lambda$ instead of all being equal to $\lambda$) and showed that the exceptional $\nu$ on a given interval $I$ with $h_\nu(\lambda)=\log M_\lambda$, understood as points in $[0, 1]^{|I|}$, lie on a lower dimensional sub-manifold (which is easy to see if all contraction parameter equal $\lambda$). Theorem \ref{main1} shows that this sub-manifold is in fact empty.

Our second result is a positive result concerning the case $\lambda$ with $|\sigma(\lambda)|<1$ for all Galois conjugates, in which case $M:=M_\lambda$ is an integer. We find that if $\nu$ satisfies certain linear equations,  maximal entropy is achieved. We show furthermore that in this case, all the nice properties we have come to expect are present. 

For $n \in \N$, let $G_n$ be the additive finite abelian group $G_n:=\Z[\beta]/\beta^n \Z[\beta]$, where $\beta=\lambda^{-1}$.   Note that the natural projection $G_n \twoheadrightarrow G_{n-1}$ induces the inclusion $\widehat{G_{n-1}} \subset \widehat{G_n}$, where $\widehat{G_n}:=\{\psi: G_n \to \mathbb{S}^1\}$ is the dual group. 
We let $P_\nu \colon \mathbb{S}^1 \to \C$, $P_\nu(z):=\EE\left[z^{\xi_0}\right]$
be the Fourier transform of $\nu$ on $\mathbb{S}^1$.  The standard Fourier transform $\widehat{\nu}$ of $\nu$ on $\R/\Z$ is related to this by
\[
P_\nu(e(\theta))=\widehat{\nu}(\theta)
\]
for any $\theta \in \R/\Z$. 

We say that the pair $\beta, \nu$ has \textit{complete vanishing at level} $m$ if for all $\psi \in \widehat{G_m}, \psi \notin \widehat{G_{m-1}}$, there exists $k=k_\psi \geq 0$ such that
\[
P_\nu(\psi(\beta^k))=0.
\]
We point out that necessarily $k_\psi<m$ because $\psi(\beta^k)=1$ for $k \geq m$. The theorem below will be reformulated later in a more general form as Theorem \ref{main2}. 
\begin{theorem} \label{main2simp}
Let $\lambda$ be an algebraic number such that for all Galois conjugates $|\sigma(\lambda)|<1$ and let $\nu$ be a finitely supported probability measure on $\Z$. 

If there is any integer $m$ with complete vanishing at level $m$, then $h_\lambda(\nu)=\log M_\lambda$, the measure $\nu_\lambda$ is absolutely continuous on $\R$ with bounded density and has power Fourier decay on $\R$.
\end{theorem}
We consider the case of $m$ being equal to one.  In this case,
\[
G_1 \cong \Z/M\Z, 
\]
where $M:=M_\lambda$,  and complete vanishing at level $1$ is achieved if and only if $\nu$ is equidistributed modulo $M$. Indeed, complete vanishing is equivalent to $\widehat{\nu}(e(\frac{k}{M}))=0$ for all $0<k<M$ with $e(x):=\exp(-2\pi i x)$, which in turn is equivalent to $\nu$ being equidistributed modulo $M$.  Thus, the $m=1$ case of Theorem \ref{main2simp} contains the aforementioned positive results about the Garsia entropy \cite{garsia1}, \cite{feng} in which $\nu=\mathrm{Unif}(\{0, \dots, M-1\})$.

We remark on the nature of the condition of complete vanishing.  For any fixed $m$,  there are $M^{m-1}(M-1)$ many $\psi$ in $\widehat{G_n} \backslash \widehat{G_{n-1}}$ (compare Lemma \ref{lem::propGn}) and for any specific $k$,  $\psi(\beta^k)$ is just some number in $\mathbb{S}^1$.  Given $m$,  the set of all probability measures $\nu$ having complete vanishing at level $m$ is thus of the form 
\[
\bigcup_{i=1}^N \{\nu: \; P_{\nu}(z)=0 \; \forall z \in E_i\}
\]
for $N=N(M_\lambda,  m)$ and for a collection of finite sets $E_i \subset \mathbb{S}^1$.  Here, the sets $E_i$ are of the form $E_i=\{\psi(\beta^{k_{i, \psi}}), \psi \in \widehat{G_m} \backslash \widehat{G_{m-1}}\}$ with $0 \leq k_{i, \psi} <m$ for all $\psi$.  As $P_\nu(z)=0$ for all $\theta \in E_i$ is a linear condition on $l^1(\Z)$,  the theorem can be understood as saying that there is a countable collection of proper subspaces in $l^1(\Z)$ such that whenever a probability measure $\nu$ lies in one of those subspaces,  $h_\lambda(\nu)$ is maximal.

When one has a criterion as in Theorem \ref{main2simp}, the first question is usually whether it is sharp. We have the following converse result.
\begin{theorem} \label{main3}
Let $\lambda$ be algebraic such that for all Galois conjugates $|\sigma(\lambda)|<1$ and let $\nu$ be a finitely supported probability measure on $\Z$. 
If $h_\lambda(\nu)=\log M_\lambda$, then there is some integer $m$ with complete vanishing at level $m$. In particular, all the other consequences of Theorem \ref{main2simp} hold.
\end{theorem}
This theorem says that any probability measure $\nu$ with maximal entropy must lie in one of the countable many subspaces corresponding to complete vanishing at some level.  Given a specific measure $\nu$ with maximal entropy,  it gives a finite list of subspaces the measure $\nu$ has to lie in.  Following the proof, this list can be determined by calculating $G_n$ and $\psi(\beta^k)$ for finitely many $n$.  The size of this list will depend on the length of the interval $\nu$ is supported on and on the biggest denominator in the finite set $\{\theta \in \Q: \widehat{\nu}(\theta)=0\}$.

This is not quite as good as one would hope to have.  Ideally, one would like this to only depend on the length of the interval $\nu$ is supported on.  Given an interval $I$ in $\Z$, one could then calculate all of the possible subspaces above and classify all probability measures with maximal entropy as a finite collection of curves in $[0,  1]^I$. 

Under the condition that two specific coefficients of the minimal polynomial of $\beta$ are coprime,  we can improve on the bounds in Theorem \ref{main3} to have this pleasant characterisation.  The bounds in this case are essentially sharp. We will see in an example below how the information from this theorem can be used to explicitly determine all measures with maximal entropy on a given interval.  
\begin{theorem}\label{easycase}
Let $\lambda$ be algebraic such that for all Galois conjugates $|\sigma(\lambda)|<1$ and let $\beta=\lambda^{-1}$.  Let $a_d x^d+ \dots+a_1 x+M$ be the minimal polynomial of $\beta$ in $\Z[X]$ and assume $\mathrm{gcd}(a_1, M)=1$. Let $\nu$ be a probability measure supported on an interval $I$ in the integers of length $|I|$ such that the support of $\nu$ is not contained in $q\Z+a$ for some $q\geq 2$ and $a \in \Z$.  Assume that $h_\lambda(\nu)=\log M_\lambda$.

If $M$ is odd,  then there is an integer $m$ with complete vanishing at level $m$ such that
\[
|I| \geq (2m-1) (M-1).
\]
If $M$ is even, the same bound holds for $m=1$; for $m \geq 2$,  it is
\[
|I| \geq (2m-2)(M-1).
\]
\end{theorem} 
When the support of $\nu$ is contained in some $q\Z+a$, one can of course replace $\xi_0 \sim \nu $ by $(\xi_0-a) q^{-1}$ and apply the theorem to this rescaled measure. 

The reason that the proof is both significantly easier and that the bounds are better is that the group $G_n$ is cyclic under the assumptions in Theorem \ref{easycase} (proved for example in Lemma \ref{lem::propGn}).  In Section \ref{sec::easycase}, we will prove Theorem \ref{easycase} and fill in the details in the example below, which demonstrates that the bounds are sharp. 
\begin{example} \label{ex::beta}
Let $\beta$ be a root of 
\[
f_\beta(x)=3x^3+4x^2+3x+5,
\]
in which case $M_\beta=5$.  
Let $E_1:=\{\frac{j}{5},  1 \leq j \leq 4\}$ and 
\[
E_2^{(1)}:=\left\{\pm \frac{1}{5} \right\} \cup \left\{\pm \frac{1+5j}{25}, 0 \leq j \leq 4 \right\}, \quad E_2^{(2)}:=\left\{\pm \frac{2}{5} \right\} \cup \left\{\pm \frac{2+5j}{25}, 0 \leq j \leq 4 \right\}.
\] 
Let $E_\nu:=\{\theta \in (0,1): \widehat{\nu}(\theta)=0\}$.  All measures $\nu$ supported on $\{0, \dots, k\}$ attaining maximal entropy can then be classified using Theorem \ref{easycase} by
\begin{itemize}
\item If $0 \leq k <4$: None.
\item If $4 \leq k<12$: Measures with $E_1 \subset E_\nu$.
\item If $12 \leq k<20$: Measures with $E_1 \subset E_\nu$ or $E_2^{(1)} \subset E_\nu$ or $E_2^{(2)} \subset E_\nu$.
\end{itemize}
At the threshold $k=12$, a concrete example is the measure $\nu$ with
 \begin{align*}
 \nu(j)&=\frac{1}{5 \varphi^2}, & j&=0, 2, 10, 12,\\
 \nu(j)&=\frac{1}{5 \varphi}, & j&=1, 5, 7, 11,\\
 \nu(j)&=\frac{1}{5}, & j&=6,
 \end{align*}
where $\varphi$ is the golden ratio.  For this measure, $E_\nu=E_2^{(2)}$,  implying maximal entropy.
\end{example}

It will become clear in the calculations how for $\beta$ as in Theorem \ref{easycase},  all measures with maximal entropy on any given interval can be classified by sets akin to $E_1$, $E_2^{(1)}$ and $E_2^{(2)}$ their Fourier transform has to vanish on.

\subsection{A more conceptual setting} \label{subsec::concept}
In this subsection, we formally introduce the more conceptual setting of the measure $\mu_\lambda$ on the contraction space $\Al$ we highlighted earlier. We also state Theorem \ref{main2}, which is a more general version of Theorem \ref{main2simp},  and the one we will prove.

To define the contraction space, we need the classification of the places and corresponding completions of a number field from algebraic number theory.  We postpone a more detailed exposition of this to Section \ref{sec::ant}. Here,  we just mention that each completion of $K=\Q(\lambda)$ comes from a \textit{place} $\nu$ of $K$, corresponding to some absolute value $|\cdot|_\nu$ on $K$.  The completion $K_\nu$ of $K$ is then the completion of $K$ with respect to the topology induced by $|\cdot|_\nu$. 

For example, for $K=\Q$,  the places are the Archimedean absolute value on $\R$ and the $p$-adic absolute values, with all completions given by $\R$ and $\Q_p$.  For more general field extensions, the completions are either $\R$, $\C$,  or finite field extensions of $\Q_p$. 
\begin{definition}[Contraction Space] \label{def::contractionspace}
Let $\lambda$ be an algebraic number.  Set $S_<:=\{\nu: |\lambda|_\nu<1\}$.  Define the contraction space $\Al$ by
\[
\Al:=\prod_{\nu \in S_<} K_\nu.
\]
We fix a Haar measure $dx$ (with respect to addition) on $\Al$ and write $\mathrm{vol}(A):=\int 1_A \; dx$ for measurable $A \subset \Al$. When we speak of absolute continuity on $\Al$, we mean with respect to $dx$.  We also fix the distance 
\[
|y|_{\Al}:=\max\left(|y_\nu|_\nu,  \; \nu \in S_<\right).
\]
The natural diagonal action $x \mapsto \lambda x:=(\lambda x_\nu)_\nu$ is a contraction on $\Al$.  By the product formula (Theorem \ref{thm::productformula}), $\mathrm{vol}(\lambda A)=M_\lambda^{-1} \mathrm{vol}(A)$.

By taking the limit with respect to the absolute value in each of the coordinates, the random variable
\[
X_\lambda:= \lim_{n \to \infty} X_n
\]
is a well-defined random variable on $\Al$.  Denote its law by $\mu_\lambda$.  
\end{definition}
We will usually denote the law of the finite random variables $\xi_j$ on $\Z$ by $\mu$ instead of $\nu$ from now on. Whether the base measure on $\Z$ is called $\mu$ or $\nu$,  $\mu_\lambda$ will always be the self-affine measure on $\Al$ and $\nu_\lambda$ will always be the self-similar measure on $\R$ as in the introduction.

As $\mu_\lambda$ projects onto $\nu_\lambda$ and the Haar measure on $\Al$ projects onto the Lebesgue measure on $\R$ in the first coordinate,  Theorem \ref{main2simp} is a corollary of the theorem we state now. 
\begin{theorem} \label{main2}
Let $\lambda$ be an algebraic number such that for all Galois conjugates $|\sigma(\lambda)|<1$ and let $\beta:=\lambda^{-1}$. Let $\mu$ be a finite measure supported on the integers. 
If there is complete vanishing at level $m$ for some integer $m$, then $h_\lambda(\mu)=\log M_\lambda$ and the measure $\mu_\lambda$ and $\mu_\beta$ are absolutely continuous with bounded density on $\Al$ and $\Ab$ respectively. 

Furthermore, on $\Ab$, there is a constant $C$ such that $\widehat{\mu_\beta}(y)=0$ for all $|y|_{\Ab} \geq C$. On $\Al$, the measure has power Fourier decay; that is, there are constants $C, \delta$ such that 
\[
|\widehat{\mu_\lambda}(y)| \leq C |y|_{\Al}^{-\delta}
\]
for any $y \in \Al$.
\end{theorem}
The power Fourier decay in Theorem \ref{main2} is significantly easier to show in the case that $\beta$ is an algebraic integer.  In the general case,  it is not conceptually harder,  but technically more involved and relies on the tools developed in Section \ref{sec::notalgint}.  We give the easier proof in the special case at the start of Section \ref{sec::fourierdecay}.

Let us briefly talk about the interplay between $\lambda$ and $\beta$. We obviously have $h_\lambda(\mu)=h_\beta(\mu)$ because the distributions of $X_n$ and
\[
\sum_{j=0}^{n-1} \xi_j \beta^j=\beta^{n-1} X_n
\]
agree; by the product formula (Theorem \ref{thm::productformula}), one also has $M_\lambda=M_\beta$. Thus, $\lambda$ has maximal entropy whenever $\beta$ has. Morally, the reason that the setting of this theorem is so nice is that the measure $\mu_\beta$ lives on a product of non-Archimedean fields, on which analysis is much more pleasant than on $\R^d$.  What we will do in the proof of Theorem \ref{main2} is to show that under the assumption of complete vanishing, $\mu_\beta$ is very well behaved on $\Ab$ and deduce everything about $\mu_\lambda$ from its connection to $\mu_\beta$.  As we see in the proof of Theorem \ref{thm::vanishcoeff}, the interplay between the Fourier transform on $\Ab$ and $\Al$ is also crucial for the entire approach of this paper. 

To conclude the introduction,  we mention the following characterisation of maximal entropy implied by our results.
\begin{corollary} \label{cor::equiv}
Let $\lambda$ with $|\lambda|<1$ be an algebraic number such that for all Galois conjugates $|\sigma(\lambda)| \neq 1$ and let $\mu$ be a finite measure supported on $\Z$.  The following are equivalent:
\begin{enumerate}
\item $h_\lambda(\mu)=\log M_\lambda$.
\item The measure $\mu_\lambda$ is absolutely continuous on $\Al$.
\item The measure $\mu_\lambda$ has power Fourier decay on $\Al$.
\item The measure $\mu_\lambda$ has any Fourier decay on $\Al$ - that is, $\widehat{\mu_\lambda}(y) \to 0$ as $|y|_{\Al}\to \infty$.
\item For all Galois conjugates,  $|\sigma(\lambda)|<1$, and a certain set of linear equations in the weights of $\mu$ are satisfied: Namely,  there is an $m$, bounded in terms of $\lambda$ and $\mu$, with complete vanishing at level $m$.
\end{enumerate}
\end{corollary} 
\begin{proof}
The implication $5. \Rightarrow 1.$ and $5. \Rightarrow 3.$ follow from Theorem \ref{main2}. The implication $1. \Rightarrow 2.$ is due to Proposition \ref{prop::ifmaxthencont}, $2. \Rightarrow 4.$ is due to the Riemann-Lebesgue lemma (Lemma \ref{lem::riemleb}) and $3. \Rightarrow 4.$ is obvious.

Lastly, the implication $4. \Rightarrow 5.$ follows from Theorem \ref{main1} and Theorem \ref{main3}: The assumption of maximal entropy in the respective proofs is only used in the form of Theorem \ref{thm::vanishcoeff}. However, to prove the result of Theorem \ref{thm::vanishcoeff}, only Fourier decay of $\mu_{\lambda}$ is needed.  
\end{proof}

\subsection{Organization of this paper}
In Section \ref{sec::ant}, we give an exposition of the results from algebraic number theory we will need throughout this paper.  All of our theorems rely on Theorem \ref{thm::vanishcoeff}, which we prove in Section \ref{sec::abscont}.  No input from other sections is needed to prove Theorem \ref{main1} and Theorem \ref{easycase}, which we do in Section \ref{sec::main1} and Section \ref{sec::easycase}, respectively. The proof of Theorem \ref{main2} relies on the approximation of the Garsia entropy at the non-Archimedean places, which we develop in Section \ref{sec::computation}.  

Section \ref{sec::notalgint} develops tools to prove Theorem \ref{main3} and power Fourier decay in Theorem \ref{main2} in the case that $\beta$ is not an algebraic integer.  In the case that $\beta$ is an algebraic integer,  they are not needed.  The proof of power Fourier decay in Theorem \ref{main2} can be found in Section \ref{sec::fourierdecay}. At the start of Section \ref{sec::fourierdecay}, we also outline the easier proof of power Fourier decay in the case that $\beta$ is an algebraic integer. The proof of Theorem \ref{main3} is given in Section \ref{sec::main3}.

\subsection{Acknowledgements}
This paper was written as part of my PhD at the University of Cambridge.  Warmest thanks go to my supervisor P\'eter Varj\'u for the good supervision,  for pushing me to write this paper in the best way it can be written, and for the many helpful comments on the drafts on the way there.  The author has received funding from the European Research Council (ERC)
under the European Union’s Horizon 2020 research and innovation programme (grant agreement No. 803711). 
\section{Places,  the Fourier transform and other notation} \label{sec::ant}
In this section,  we recall the concepts from algebraic number theory we need throughout the paper.  We also pin down some notation.  All of the material in this section is standard and can be found for example in \cite{fouriernf}.

Throughout this paper, we are considering a number field $K$ and are interested in its completions.   Each completion $K_\nu$ corresponds to an absolute value $|\cdot|_\nu$,  as we outline below.  The equivalence class of all absolute values corresponding to the same completion is then called a \textit{place} $\nu$ of $K$,  so places and completions are in 1-1 correspondence. 

To classify all places,  let an algebraic $\lambda$ of degree $d$ be given and set $K:=\Q(\lambda)$.  On the one hand, there are the Archimedean completions of $K$ coming from the field embeddings $\sigma: K \to \R, \C$. If $\sigma$ is real,  the corresponding completion is $\R$ with the usual absolute value.  When $\sigma$ is complex-valued,  the completion is $\C$ and the topology is induced by 
$|x|_\nu:=| \sigma(x) |^2_\C$ (one of each pair of conjugates is chosen).  This is technically not an absolute value because it fails the triangle inequality, but this will not matter much.  This normalisation is customary in algebraic number theory to avoid having an exponent of two in the product formula below. 

On the other hand, there are the non-Archimedean absolute values that correspond to a prime ideal of the ring of integers $\OK$. For $K=\Q$, these are all the prime numbers in $\Z$ and the corresponding completions are the $p$-adic fields $\Q_p$. In this case, writing a rational number as $x=\frac{a}{b} p^s$, $a, b$ coprime to $p$, the corresponding absolute value is given by $|x|_p:=p^{-s}$.

For more general $K$, any prime ideal $P$ in $\OK$ lies over a prime $p \in \Z$. For $\nu$ corresponding to $P$, the absolute value $|\cdot|_\nu$ is an extension of $|\cdot|_p$ and the completion $K_\nu$ is a finite field extension of $\Q_p$.  For non-Archimedean $\nu$ one has the ultrametric inequality
\[
|x+y|_\nu \leq \max(|x|_\nu, |y|_\nu).
\]
We omit to talk about the normalisation of the absolute values here - they are normalised such that the product formula below holds. If the reader is interested in the normalisation,  details can be found in Chapter 5 of \cite{fouriernf}.  

The collection of the absolute values above then consists of all places of $K$ (this is Ostrovski's  Theorem).   Another important result is the product formula,  which has frequently been a key input of proofs related to the Garsia entropy (see for example \cite{garsia1}, \cite{feng}, \cite{bv}).
\begin{theorem}[Product formula; \cite{fouriernf}, Theorem 5-14 (i)] \label{thm::productformula}
For any non-zero $x \in K$,  
\[
\prod_{\nu} |x|_\nu=1,
\]
where $\nu$ runs over all places of $K$.
\end{theorem}
The product is well defined because,  for each non-zero $x \in K$, there are at most finitely many $\nu$ with $|x|_\nu \neq 1$.  We now define the Fourier transform on the completions $K_\nu$.  

Fix a field extension $K_\nu$ of $\Q_p$. One can consider the additive group of $K_\nu$ as a locally compact Abelian group with characters and a dual group. It was shown in Tate's thesis that the fields $\Q_p$ and $K_\nu$ are self-dual \cite{tate}. Together with the Haar measure on $K_\nu$, the characters on $K_\nu$ induce an isomorphism of $L^2(K_\nu)$ onto itself in the usual way. One defines the Fourier transform of $f$ by fixing an isomorphism of $\widehat{K_\nu}$ with $K_\nu$ via $\widehat{K_\nu}=\{\psi_y, y \in K_\nu\}$ and setting
\[
\widehat{f}(y):=\int_{K_\nu} f(x) \psi_y(x) dx,
\]
where $dx$ is the Haar measure. Roughly speaking, the Fourier transform on $K_\nu$ is like the Fourier transform on $\R$, only nicer.

For concreteness, we also detail how this identification looks like: The standard character of $\Q_p$ is given by 
\[
\psi^{\Q_p}(y):=e\left(-\sum_{i=-N}^{-1} a_i p^i \right),
\]
for any element $y=\sum_{i=-N}^\infty a_i p^i \in \Q_p$, $a_i \in \{0, \dots, p-1\}$, where $e(x):=\exp(-2 \pi i x)$.
For the identification of $\widehat{\Q_p}$ with $\Q_p$ one sets 
\[
\psi^{\Q_p}_y: \Q_p \to \mathbb{S}^1, \quad \psi^{\Q_p}_y(x)=\psi^{\Q_p}(xy).
\]
On the finite extension $K_\nu$ of $\Q_p$, one defines the standard character via the field trace by
\[
\psi^{K_\nu}(x):=\psi^{\Q_p}(\tr_{K_\nu/\Q_p}(x))
\]
and $\psi^{K_\nu}_y: K_\nu \to \mathbb{S}^1$ via $\psi^{K_\nu}_y(x):=\psi^{K_\nu}(xy)$. We remark that $\psi^{K_\nu}(x)=1$ if $|x|_\nu \leq 1$.  For the Archimedean fields, the standard character is given by $\psi^\R(x)=e(x)$ and $\psi^\C(x)=e(2\mathrm{Re}(x))$.  

We set 
\[
S_<:=\{\nu: \; |\lambda|_\nu<1\}, \quad \Si:=\{\nu \in S_<: \; \nu \mathrm{\; arch.}\}, \quad  \Sf:=\{\nu \in S_<: \; \nu \mathrm{\; non-arch.}\} 
\]
and also $S_>:=\{\nu: \; |\lambda|_\nu>1\}$,  $S_=:=\{\nu: \; |\lambda|_\nu=1\}$.

We recall that the contraction space is defined by $\Al:=\prod_{\nu \in S_<} K_\nu$ and define analogously
\[
\Ali:=\prod_{\nu \in \Si} K_\nu, \quad \Alf:=\prod_{\nu \in \Sf} K_\nu.
\]
We note that $\Ali \cong \R^s$ for some $s$ and that $\Ab=\prod_{\nu \in S_>} K_\nu$, where as always $\beta=\lambda^{-1}$.  

There is a natural diagonal embedding $D: K \to \Al$.  We will frequently drop the $D$ and just write $y \in \Al$ in place of $D(y) \in \Al$.  

The standard character $\psi^{\Al}: \Al \to \mathbb{S}^1$ on $\Al$ is given by 
\[
\psi^{\Al}(x):=\prod_{\nu \in S_<} \psi^{K_\nu}(x_\nu)
\]
for $x=(x_\nu)_{\nu \in S_<} \in \Al$. Once again, this induces an isomorphism of $\widehat{\Al}$ to $\Al$ via
\[ 
\psi^{\Al}_y: \Al \to \mathbb{S}^1, \quad \psi^{\Al}_y(x)=\psi^{\Al}(xy)
\]
for $y \in \Al$.  One defines the standard character and duality of $\Alf$ and $\Ali$ in the same way,  with the product ranging over $\Sf$ and $S_>$ instead.

The Fourier transform on $\Al$ is defined by setting 
\[
\widehat{f}(y):=\int_{\Al} f(x) \psi^{\Al}_y(x) dx.
\]
One has the direct analogue of the Riemann-Lebesgue lemma, which one can prove in exactly the same way as the classical Riemann-Lebesgue lemma, exploiting that simple functions are dense in $L^1(\Al)$. 
\begin{lemma}[Riemann-Lebesgue] \label{lem::riemleb}
For any $f \in L^1(\Al)$, $\widehat{f}(y) \to 0$ as $y \to \infty$. 
\end{lemma} 
An important input in the proofs of the lack of absolute continuity for Pisot numbers (\cite{garsia1}, \cite{erdos}) and algebraic units (\cite{bv}) for the Bernoulli measure is the fact that for $x$ an algebraic integer,  
\[
\prod_{\nu \; \mathrm{arch}} \psi^{K_\nu}(x)=e(\tr_{K/\Q}(x))=1.
\]
This is not true for $x$ which is not an algebraic integer,  but the proposition below provides an analogous statement we will use frequently. Because $\psi^{K_\nu}(x)=1$ if $|x|_\nu \leq 1$, the product is well-defined.
\begin{prop}[\cite{fouriernf}, Proposition 7-15 (ii)] \label{prop::adelefact}
For any $x \in K$, 
\[
\prod_{\mathrm{all} \; \nu} \psi^{K_\nu}(x)=1,
\]
where $\nu$ runs over all places of $K$.
\end{prop} 
We also let for $y \in K$
\[
\Sigma_<(y):=\left(\sum_{\sigma \in \Si \; \mathrm{real}} \sigma(y)+\sum_{\sigma \in \Si \; \mathrm{complex}} 2 \mathrm{Re}(\sigma(y)) \right).
\]

We can identify $\Ali$ with some $\R^s$ by identifying $\C$ with $\R^2$ via $z \mapsto(\mathrm{Re}(z),  \mathrm{Im}(z))$.  We write $D_\infty: K \to \R^s$ for the map consisting of the diagonal embedding $K \to \Ali$ composed with this identification.  We will also frequently drop the $D_\infty$ and just write $y \in \R^s$ in place of $D_\infty(y) \in \R^s$.  

This then naturally induces a bilinear form on $\Ali$ given by 
\[
\langle x, y \rangle_<=\Sigma_<(xy)=\sum_{\nu \; \mathrm{real}} x_\nu y_\nu+2\sum_{\nu \; \mathrm{complex}} \mathrm{Re}(x_\nu y_\nu)
\]
which makes $\psi_x^{\Ali}(y)=e(\langle x, y \rangle_<)$.  

This is almost the standard euclidean scalar product on $\R^s$, but not quite. We let $S:\R^s \to \R^s$ be the map corresponding to $z \mapsto 2 \overline{z}$ in all complex coordinates on $\Ali$ - explicitly,  $e_i \mapsto e_i$ if $i$ corresponds to a real embedding,  $(e_i, e_{i+1}) \mapsto 2 (e_i, -e_{i+1})$ if $i,  i+1$ corresponds to a complex embedding. Then,
\[
\langle x, y \rangle_<=\langle S(x),  y \rangle_{\mathrm{std}},
\]
where $\langle - , -\rangle_{\mathrm{std}}$ is the standard scalar product on $\R^s$.  
Because we only use $\langle -, -\rangle_<$ to induce the duality of $\R^s$, the lack of positive definiteness of $\langle - , -\rangle_<$ will be immaterial to us.  We then have for all $x, y \in K$ that
\[
\psi^{\Al}_x(y)
=e(\Sigma_<(xy))\psi^{\Alf}_x(y)=e(\langle x,  y \rangle_<)\psi^{\Alf}_x(y).
\]

To finish, we define the polar lattice.  A lattice $L$ in $\R^s$ is a discrete subgroup of full rank.  We define the \textit{polar lattice} to be
\[
L^*:=\{x \in \R^s: \; \langle x, y \rangle_< \in \Z \; \forall y \in L \}.
\]

There are further tools from algebraic number theory we need in Section \ref{sec::notalgint} and after,  which will be introduced then.

Lastly, we pin down the big-$O$ notation: For some domain $U$ and functions $f, g \colon U \to \R$,  we write $f=O(g)$ if there is a constant $C$ such that $|f(u)| \leq C |g(u)|$ for all $u \in U$.  We write $O_S$ if the constant $C$ is allowed to depend on some quantity $S$.

\section{The general approach} \label{sec::abscont}
In this section,  we outline the general approach,  which is used as a starting point in all of our theorems.  The main goal of this Section is to prove Theorem \ref{thm::vanishcoeff} and
Proposition \ref{prop::ifmaxthencont}.  Both are more general than what appeared in the literature before (compare \cite{garsia1} and \cite{bv}, Proposition 17).  But proving them is the classical approach to showing absolute continuity or the lack thereof for the Bernoulli measure (\cite{garsia1}, \cite{garsia2}, \cite{bv}, \cite{feng}) and the ideas ultimately go back to Garsia and Erd\H os.

We set
\[
\CC:=\left\{x \in K \; \Big| \; \forall s \in \Z: \psi^{\Al}_x(\lambda^s)=\psi^{\Ab}_x(\lambda^s)^{-1} \right\}=\left\{x \in K \; \Big| \; \forall s \in \Z: \prod_{\nu \in S_=} \psi^{K_\nu}_x(\lambda^s)=1  \right\},
\]
where as always $\beta=\lambda^{-1}$ and $S_==\{\nu: |\lambda|_\nu=1\}$.  The second equality is due to 
Proposition \ref{prop::adelefact}.  We recall that $P_\mu(e(\theta))=\widehat{\mu}(\theta)$ for $\theta \in \R/\Z$.

\begin{theorem} \label{thm::vanishcoeff}
Let $\lambda$ be algebraic without a conjugate on the unit circle and let $\mu$ be a measure on the integers such that $h_\lambda(\mu)=\log M_\lambda$.  Then $\mu_\lambda$ and $\mu_\beta$ and are absolutely continuous with bounded density on $\Al$ and $\Ab$, respectively. Moreover,  for any $x \in \CC$ there is $k=k(x) \in \Z$ such that
\[
P_\mu(\psi^{\Al}_x(\lambda^k))=0.
\]
\end{theorem}
Having a result like this is not of much use unless we know that $\CC$ contains non-zero elements. 
If $\lambda$ is a Salem number,  there is by definition an Archimedean place $\nu_s$ in $S_=$.  This renders this approach pretty much useless because there are very few $x$ for which $\prod_{ \nu \in S_=} \psi^{\Al}_x(\lambda^s)=1$, even for fixed $s$.  This keeps Salem numbers out of reach for the methods in this paper.  As proving pretty much any properties of $\nu_\lambda$ is a wide open problem for Salem numbers, this may not be too much of a surprise. 

When $\lambda$ is not a Salem number and there are only non-Archimedean places in $S_=$, the situation is much better. Recalling that $\psi^{K_\nu}(y)=1$ if $|y|_\nu \leq 1$ for non-Archimedean $\nu$,  one sees that for example $\Z[\lambda,  \lambda^{-1}] \subset \CC$. 

In the papers \cite{erdos}, \cite{garsia1}, \cite{bv}, where one is dealing with the Bernoulli measure $\nu^{\mathrm{Ber}}$ and $\lambda$ a Pisot number or an algebraic unit, one is done at this stage: The theorem can be applied for $x=1$. Indeed, for say a Pisot number,
\[
\psi^{\Al}_x(\lambda^k)=e(\lambda^k) \neq -1
\]
for any $k$, while the only zero of $P_{\nu^{\mathrm{Ber}}}$ is at $-1$. 

That $x=1$ works is quite specific to the Bernoulli measure though.   Given a general measure $\nu$,  $P_\nu$ may well have a zero at some $e(\lambda^k)$. To deal with a general $\nu$, we will need to tap into the full breadth of the $x$ available to us from the corollary. The approach of using all of the vanishing Fourier coefficients instead of only $x=1$ is the main difference between our work and all of the prior work quoted above.

We now sketch how this theorem is proven. The most important building block of the proof is the proposition below.
\begin{prop} \label{prop::ifmaxthencont}
Let $\lambda$ be algebraic without a Galois conjugate on the unit circle. If $h_\lambda(\mu)=\log M_\lambda$, then $X_\lambda$ is absolutely continuous on $\Al$ with bounded density.  (The converse also holds - see Corollary \ref{cor::equiv}). 
\end{prop}

The second ingredient is the Riemann-Lebesgue Lemma (Lemma \ref{lem::riemleb}), which implies that if $\mu_\lambda$ is absolutely continuous, then $\widehat{\mu_\lambda}(y) \to 0$ as $y \to \infty$ in $\Al$. 

One can explicitly express the Fourier transform as an infinite product.  One then considers for $x$ given the elements $\lambda^{-n} x$, which converge to infinity in $\Al$.  All values of $\widehat{\mu_\lambda}(\lambda^{-N}x)$ will be bounded above by a bi-infinite product. The key is then that for $x \in \CC$, both sides of this bi-infinite product converge. This then implies that there must be a vanishing coefficient somewhere in the product, finishing the sketch of the proof.\\

For the rest of the section,  we prove the results above.  We start by proving a lemma we need in the proof.  We then prove Proposition \ref{prop::ifmaxthencont}, which is the main part of this section.  We end this section by giving a formal proof of Theorem \ref{thm::vanishcoeff}.

The main idea in this section and in Section \ref{sec::computation} will be to 'put a box on it', that is, to smoothen out the self-similar measure by fixing some scale and identifying points which are close together (whether in the Archimedean or in the non-Archimedean setting). We start with a preparatory lemma also used in Section \ref{sec::computation}.  This lemma also implies that $h_\lambda(\mu) \leq \log M_\lambda$. 

Recall
\[
X_n=\sum_{j=0}^{n-1} \lambda^j \xi_j \in \Al,
\]
where $X_n \in K$ is embedded diagonally and $\xi_j \sim \mu$.
\begin{lemma} \label{lem::separation}
Let $\lambda$ be an algebraic integer without a Galois conjugate on the unit circle. There is an absolute constant $\eta=\eta(\lambda, \mu)>0$ such that if $x, y \in \mathrm{supp}(X_n)$ are such that 
\[
|x-y|_{\nu} \leq 2 \eta |\lambda|_\nu^n
\]
for all $\nu \in S_<$, then $x=y$. In particular, there exists an absolute constant $C=C(\lambda, \mu)$ such that 
\[
|\mathrm{supp}(X_n)| \leq C M_\lambda^n.
\]
\end{lemma}
\begin{proof}
Recall $\beta=\lambda^{-1}$ and
\[
S_>=\{\nu: |\lambda|_\nu>1\}=\{\nu: |\beta|_\nu<1\}.
\]
Let $z \in \mathrm{supp}(X_n)-\mathrm{supp}(X_n)$ and let $\mu$ be supported on $\{0, \dots, B\}$. Note that then 
\[
\beta^{n-1} z \in \{\sum_{n=0}^\infty a_n \beta^n, 0 \leq a_n \leq B\} \subset \Ab.
\]
As
\[
\left|\sum_{n=0}^\infty a_n \beta^n \right|_\nu \leq \max_{0 \leq j \leq B}\{ |j|_\nu\} \sum_{n=0}^\infty |\beta|_\nu^n=\frac{\max_{0 \leq j \leq B}\{ |j|_\nu\}}{1-|\beta|_\nu}
\]
for all $\nu \in S_>$, we see that there is an absolute $C_0=C_0(\lambda, \mu)$ such that 
\[
|z|_\nu=|\lambda|_\nu^{n-1} |\beta^{n-1} z|_\nu \leq C_0 |\lambda|_\nu^{n-1}
\]
for all $\nu \in S_>$. For any $\nu$ with $|\lambda|_\nu=1$, which by assumption must be non-Archimedean, $|z|_\nu \leq 1$ by the strong triangle inequality. 

Then, by the product formula (Theorem \ref{thm::productformula}), for any distinct $x, y$ as in the lemma,
\[
2 \eta M_\lambda^{-n} \geq \prod_{\nu \in S_<} |x-y|_\nu=\prod_{\nu \notin S_<} |x-y|^{-1}_\nu \geq \prod_{\nu \in S_>} |x-y|^{-1}_\nu \geq \prod_{\nu \in S_>} C_0^{-1} |\lambda|^{-(n-1)}_\nu=C_0^{s} M_\lambda^{-(n-1)},
\]
for $s=|S_>|<\infty$. This is impossible for $\eta$ small enough.

For the second part, we let 
\begin{equation} \label{eq::boxdef}
Q_r=\prod_{\nu \in S_<} \{x: |x|_\nu \leq r\} \subset \Al 
\end{equation}
be the box with side length $r$ and note that
\[
\lambda Q_r=\prod_{\nu \in S_<} \{x: |x|_\nu \leq r|\lambda|_\nu\}.
\]
We know that $X_\lambda $ is contained in $Q_N$ for some big absolute $N$. By the first part of the lemma, $x+\lambda^{n-1} Q_{\eta}$ and $y+\lambda^{n-1} Q_{\eta}$ are disjoint for all distinct $x, y \in \mathrm{supp}(X_n)$. Thus,
\[
|\mathrm{supp}(X_n)| \leq \frac{dx(Q_N)}{dx(\lambda^{n-1} Q_{\eta})}=\frac{dx(Q_N)}{M_\lambda^{-(n-1)} dx(Q_{\eta})},
\]
which shows the claim.
\end{proof}
We now show Proposition \ref{prop::ifmaxthencont}, which asserts that under the assumption of maximal entropy, $\mu_\lambda$ is absolutely continuous on $\Al$.  The proof is almost identical to the proof of Proposition 17 in \cite{bv},  with the ideas going back to Garsia.
\begin{proof}[Proof of Proposition \ref{prop::ifmaxthencont}]
It is a consequence of the Jessen-Wintner Theorem (Theorem 35 in \cite{jessenwintner}) that $\mu_\lambda$ is either singular or absolutely continuous with respect to the Haar measure.  Assume for contradiction that $X_\lambda$ is singular. There is an absolute $C_0=C_0(\lambda, \mu)$ such that
\[
|X_\lambda-X_n|_\nu \leq \max_{\omega}|\xi_0(\omega)|_\nu \sum_{k=n}^\infty |\lambda|^k_\nu \leq C_0 |\lambda|_\nu^n
\]
for all $\nu \in S_<$. 

Fix $\eps>0$. As $X_\lambda$ is singular with respect to the Haar measure $dx$, there is a closed set $B_0 \subset \Al$ with $dx(B_0)<\eps$ such that $\PP(X_\lambda \in B_0) \geq 1-\eps$. Note that as $B_0$ is closed, for $n^\prime$ big enough, $dx(B_1)\leq 2\eps$, where
\[
B_1:=B_0+ 2\lambda^{n^\prime} Q
\]
with $Q:=Q_1$ as in (\ref{eq::boxdef}).  Set 
\[
B:=B_0+\lambda^{n^\prime} Q
\]
and let $n$ be so big that $C_0 |\lambda|_\nu^n \leq \eta | \lambda|_\nu^{n^\prime}$ for all $\nu \in S_<$, where $\eta$ is as in Lemma \ref{lem::separation}. We then see that $\PP(X_n \in B)\geq \PP(X_\lambda \in B_0) \geq 1- \eps$ by the estimate.

On the other hand, the sets $x+\lambda^n Q_\eta$ are disjoint for all $x \in \mathrm{supp}(X_n)$ by Lemma \ref{lem::separation}. As also by construction, 
\[
B+\lambda^n Q_\eta \subset B+ \lambda^{n^\prime} Q \subset B_1, 
\]
we deduce that
\[
2 \eps \geq dx(B_1) \geq dx(\lambda^n Q_\eta) |\mathrm{supp}(X_n) \cap B|=dx(Q_\eta) M_\lambda^{-n} |\mathrm{supp}(X_n) \cap B|.
\]
Thus, $|\mathrm{supp}(X_n) \cap B| \leq C_1 \eps M_\lambda^n$ for some absolute $C_1$. 

Recall that $|\mathrm{supp}(X_n)| \leq CM_\lambda^n$ by Lemma \ref{lem::separation} and that $\PP(X_n \notin B) \leq \eps$.  Set $Y:=\mathbf{1}_{X \in B}$ and note that $H(Y) \leq \log 2$.
\begin{align*}
H(X_n)&=\PP(X_n \in B) H(X_n| X_n \in B)+\PP(X_n \notin B) H(X_n|X_n \notin B)+H(Y) \\
&\leq \PP(X_n \in B) \log C_1 \eps M_\lambda^n+\PP(X_n \notin B)\log C M_\lambda^n+\log 2\\
&\leq n \log M_\lambda +\log C_1 C+(1-\eps)\log \eps+\log 2 .
\end{align*}
The expression on the right-hand side is strictly smaller than $n \log M_\lambda$ for $\eps$ small enough.  It is well-known and easy to see that the sequence $H(X_m)$ is subadditive, which means that $\frac{1}{m} H(X_m) \to h_\lambda(\mu)$ from above.  Thus, the above implies that
\[
h_\lambda(\mu) \leq \frac{1}{n}H(X_n)<\log M_\lambda,
\]
contradiction.
\end{proof}
We now conclude this section by proving the theorem.
\begin{proof}[Proof of Theorem \ref{thm::vanishcoeff}]
As we noted below Theorem \ref{main2},  $h_\lambda(\mu)=h_\beta(\mu)$ and $M_\lambda=M_\beta$, so by Proposition \ref{prop::ifmaxthencont},  $\mu_\lambda$ and $\mu_\beta$ are absolutely continuous on $\Al$ and $\Ab$, respectively. 

Applying  the Riemann-Lebesgue lemma (Lemma \ref{lem::riemleb}), we deduce that $\widehat{\mu_\lambda}(y) \to 0$ for any $y \to \infty$.  The Fourier transform of $\mu_\lambda$ is 
\begin{align*}
\widehat{\mu_\lambda}(x)=\EE\left[\psi^{\Al}_x\left(\sum_{j=0}^\infty \xi_j \lambda^j \right) \right]=\prod_{j=0}^\infty P_\mu(\psi^{\Al}_x(\lambda^j)).
\end{align*}
We pick a non-zero $x \in \Al$ and consider the sequence $\lambda^{-n}x \to \infty$.  Recalling that $\psi^{\Al}_{xy}(z)=\psi^{\Al}_{x}(yz)$ by definition,  we see that
\[
0=\lim_{n \to \infty} |\widehat{\mu_\lambda}(\lambda^{-n}x)|=\lim_{n \to \infty} \prod_{j=-n}^\infty |P_\mu(\psi^{\Al}_x(\lambda^j))|= \prod_{j=-\infty}^\infty |P_\mu(\psi^{\Al}_x(\lambda^j))|.
\]
We now show that for $x \in \CC$, both tails of this infinite product converge. This concludes the proof of the theorem because the only way that the product can vanish is then if one of the factors vanishes.

The convergence of the positive tail is true for any $x$ and should not be very surprising: After all, the random variable $X_\lambda$ converges in $\Al$.  Explicitly,  one sees this because on the one hand,  for non-Archimedean $\nu \in S_<$,  $\psi^{K_\nu}(\lambda^n x)=1$ for $n=n(x)$ sufficiently big (because then $|\lambda^n x|_{\nu} \leq 1$).  On the other hand,  the Archimedean contribution $\prod_{n=0}^\infty P_\mu(e(\langle x, \lambda^n \rangle_<))$ can be bounded because
\[
\sum_{n=0}^\infty |\langle x,  \lambda^n \rangle_<| \leq 2 \sum_{\sigma \in \Si} |\sigma(x)| \sum_{n=0}^\infty  |\sigma(\lambda)|^n< \infty.
\]

The key point is now that for $x \in \CC$, the negative tail
\[
\prod_{j=-\infty}^{-N}P_\mu(\psi^{\Al}_x(\lambda^j))=\prod_{j=N}^{\infty} P_\mu(\psi^{\Al}_x(\beta^{j}))=\prod_{j=N}^{\infty} P_\mu(\psi^{\Ab}_x(-\beta^{j}))
\]
can be bounded in exactly the same manner as the positive tail.  
\end{proof}

\section{Proof of Theorem \ref{main1}} \label{sec::main1}
We start by proving Theorem \ref{main1} in the case of a Pisot number because the main idea comes across more clearly.  We then give the proof of the general case.
\begin{proof}[Proof of Theorem \ref{main1}, Pisot case] \label{sketch1}
Assume that $\mu$ is such that $h_\lambda(\mu)=\log M_\lambda$.  We assume that $\lambda$ is the inverse of a Pisot number, in which case $\Al=\R$.  We also assume that it is an algebraic unit.  Let $E:=\{\theta \in \R/\Z: \widehat{\mu}(\theta)=0\}$, which is finite because the Fourier transform of $\mu$ is a polynomial in $e(\theta)$. By Theorem \ref{thm::vanishcoeff}, we know that for any $x \in \Z[\lambda]$ there is a $k \in \Z$ such that 

\[
\widehat{\mu}(x \lambda^k)=P_\mu(\psi^{\Al}_x(\lambda^k))=0.
\]

This condition is equivalent to saying that for each $x \in \Z[\lambda]$ there is $k \in \Z$ such that $x \lambda^k \in E+\Z$. Any $x \lambda^k$ can be written as $\sum_{j=0}^{d-1} b_j \lambda^j$ for some $b_j \in \Z$ using the minimal polynomial of $\lambda$. Because $\lambda$ has degree $d$, all of those expressions are distinct, which in turn implies that all $\sum_{j=1}^{d-1} b_j \lambda^j$ are distinct modulo $\Z$. As $E$ is finite and $0 \notin E$, there is some big number $N$ such that $\sum_{j=1}^{d-1} N b_j \lambda^j \notin E+\Z$ for any choice of $b_j$.  If we set $x=N$, we can thus ensure that for any value of $k$, $x \lambda^k \notin E+\Z$, contradiction.
\end{proof}

In the general case, we proceed similarly: Theorem \ref{thm::vanishcoeff} gives us for any $x \in \Z[\lambda]$ a $k \in \Z$ such that
\[
0=P_\mu(\psi^{\Al}_x(\lambda^k))=P_\mu\left(e(\Sigma_<(x\lambda^k)) \psi^{\Alf}_x(\lambda^k)\right).
\]
Here,  as defined in Section \ref{sec::ant},  $\psi^{\Alf}_x=\prod_{\Sf} \psi_x^{K_\nu}$ is the contribution of the non-Archimedean completions and $\Sigma_<(xy)$ is an argument of the Archimedean contributions $\psi^{\Ali}_x(y)$.  One has that $\psi^{\Alf}(x \lambda^k)$ is of the form $e(q)$ for some $q \in \Q$ for any $x$ and $k$, which makes its contribution quite unimportant for this argument.  By a similar, if technically more involved argument,  as in the case of a Pisot number,  one can prove the lemma below.  We will apply the lemma below with $E$ taken to be the zeroes of $\widehat{\mu}$.

\begin{lemma} \label{lem::irrational}
Let $\lambda$ be algebraic without Galois conjugates on the unit circle.  Let $E \subset \R/\Z$ be a finite set with $0 \notin E$. Then, there exist a prime $p$ depending only on $\lambda$ and $E$ such that if $w \in p \Z[\lambda]$ is such that $\Sigma_<(w \lambda^k) \notin \Q$, then $\Sigma_<(w \lambda^k) \notin E+\Q$ (irrespective of k).
\end{lemma}
\begin{proof}[Proof of Lemma \ref{lem::irrational}]
We define $u_j:=\Sigma_<(\lambda^j)$ for $0 \leq j < d$ and let $w_0:=u_0, \dots, w_r$ be a maximal $\Q$-linearly independent subset of the $u_j$. Let $W:=\langle w_j, j \leq r\rangle_\Q \subset \R$ be the space generated by them. Because $1, \dots, \lambda^{d-1}$ span $K$, one has that $\Sigma_<(K)=W$. Any elements $x$ of $E$ not contained in $W$ cannot contain any elements of the form $\Sigma_<(w \lambda^k)+\Q$, so without loss of generality,  there are none. 

Take $z \in E$, which can uniquely be written as
\[
z=\sum_{j=0}^r z_j w_j, \quad z_j \in \Q.
\]
Remembering that $w_0=|\Si| \in \Q$, we see that 
\[
z+\Q \subset \Q+\sum_{j=1}^r z_j w_j.
\]
The union of the respective coefficients $z_j \in \Q$ for all $z \in E$ is finite, so we can find a big prime $p$ which appears in neither the numerator nor the denominator of any $z_j$. This choice of $p$ implies that any element in $W$ which has a coefficient with $p$ in the numerator in front of one of $w_1, \dots, w_r$ will not be in $E+\Q$.

Apart from this,  we can also impose on $p$ that it is coprime to the first and last coefficient in the minimal polynomial of $\lambda$ and that it does not appear in any of the denominators when we express the $u_j$ as linear combinations of the $w_j$ over $\Q$. 

The prime $p$ chosen like that will be the one with which the conclusion holds: Take $w$ as in the statement
and write $\lambda^k w=\sum_{j=0}^{s} p c_j \lambda^{k+j}, c_j \in \Z$. We can use the minimal polynomial of $\lambda$ to write 
\[
\lambda^{d}=-\frac{1}{a_d}(\sum_{i=0}^{d-1} a_i \lambda^i) \quad \mathrm{and} \quad
\lambda^{-1}=-\frac{1}{a_0}(\sum_{i=1}^{d} a_{i} \lambda^{i-1}).
\]
Using this iteratively, we see that
\[
\lambda^k w=\sum_{j=0}^{d-1} p d_j \lambda^j
\]
where the $d_j$ are in $\Q$ with denominators dividing $a_0^N a_d^N$ for some big $N$. 

By our choice of $p$, none of them can cancel $p$ out. Thus,
\[
\Sigma_<(w \lambda^k)=\sum_{j=0}^{r} p e_j w_j
\]
for some $e_j \in \Q$ which, by our choice of $p$, have denominators coprime to $p$. By the property above, $\Sigma_<(w \lambda^k) \in E+\Q$ can only happen if $e_1=\dots=e_r=0$. This is exactly equivalent to $\Sigma_<(w \lambda^k) \in \Q$.
\end{proof}

We sketch how one concludes the proof from here.  If $S_<$ were to consist of all or no embeddings, this lemma does not help us in any way in the proof of Theorem \ref{main1} because $\Sigma_< =\mathrm{tr}_{K/\Q}$ or $\Sigma_<=0$, which is $\Q$-valued (and indeed Theorem \ref{main1} is then false).  However, as the lemma below shows,  $\Sigma_<$ being $\Q$-valued happens only in these cases.
\begin{lemma} \label{lem::obs}
Let $S$ be a set of field embeddings of $\lambda$ into $\C$ (picking both or none of a pair of complex embeddings). If $\sum_{\sigma \in S} \sigma(\lambda^j)\in \Q$ for $d$ many consecutive $j$, then $S$ must consist of all or no field embeddings.
\end{lemma}
\begin{proof}[Proof of Lemma \ref{lem::obs}]
(The neat proof of this lemma is due to Wojtek Wawrów, see  \cite{wawrow}) Assume that the assumption is satisfied. Using the minimal polynomial of $\lambda$, we see that $\sum_{\sigma \in S} \sigma(\lambda^j) \in \Q$ for all $j \in \Z$. Consider the polynomial
\[
g(x):=\prod_{\sigma \in S} (x-\sigma(\lambda)),
\] 
which has coefficients of elemantary polynomials in $\sigma(\lambda)$ for $\sigma \in S$. The first $d$ elementary polynomials can be expressed in terms of the power polynomials $\sum_{\sigma \in S} \sigma(\lambda)^j$ for $j \leq d$ by Newton's identities. So under our assumption on the power polynomials, one has $g \in \Q[x]$.  As $g(\lambda)=0$, $g$ must then be divisible by the minimal polynomial of $\lambda$.  So either $g=0$ (no field embeddings) or $\mathrm{deg}(g) \geq d$ (all field embeddings). 
\end{proof}
We now conclude the proof of Theorem \ref{main1}.
\begin{proof}[Proof of Theorem \ref{main1}] 
We let $E \subset (0,1)$ be the zeroes of $\widehat{\mu}$, which are finitely many because $P_\mu$ is a polynomial and $e(E)$ are the zeroes of $P_\mu$. 
We recall that by Theorem \ref{thm::vanishcoeff},  under the assumption of maximal entropy,  for any $v \in \Z[\lambda]$ there is a $k \in \Z$ such that $P_\mu(\psi^{\Al}_v(\lambda^k))=0$.  Because $\psi^{\Alf}$ takes values in $e(\Q)$, this implies that $\Sigma_<(v\lambda^k) \in E+\Q$.  

We fix $p$ as in Lemma \ref{lem::irrational}. The claim below shows that one can find a reasonably small $v$ such that for a large stretch of $k$,  Lemma \ref{lem::irrational} prevents the vanishing coefficient to occur.  Once the claim is shown, the rest of the proof will be straight forward estimate of the size of $\Sigma_<(v\lambda^k)$ for large $k$. 
\begin{claim} \label{claim::main1help}
There is a constant $c_0$ depending only on $\lambda$ such that for any $K$ sufficiently big,  there is a $v \in p\Z[\lambda]$ of the form $v=\sum_{i=0}^{d-1} p v_i \lambda^i$ with $v_i \in \Z$, $|v_i| \leq K$ such that $\Sigma_<(v \lambda^k) \notin \Q$ for $|k| \leq c_0 K$.
\end{claim}
\begin{proof}[Proof of Claim \ref{claim::main1help}]
Fix an integer $k$ and set 
\[
\varphi_k: \Q^d \to W, \quad \vec{v}=(v_0, \dots, v_{d-1}) \mapsto \sum_{j=0}^{d-1} p v_j \Sigma_<\left(\lambda^{j+k} \right).
\]
By our assumption in Theorem \ref{main1}, $\Si$ does not consist of all or no field embeddings. Thus,  by Lemma \ref{lem::obs}, $\varphi_k^{-1}(\Q) \neq \Q^d$.  Set $A_k:=\varphi_k^{-1}(\Q) \cap \Z^d$. It now just remains to count points in $\Z^d$.

The ball $B_K:=\{\vec{v} \in \Z^d\big| \forall j: |v_j| \leq K\}$ contains at least $K^d$ many points. For any $A \subset \Z^d$ that is contained in a proper subspace of $\Q^d$,
\[
|A \cap B_K| \leq C K^{d-1},
\]
where $C$ is some absolute constant depending only on the dimension $d$. Thus, by the union bound, 
\[
\left|\bigcup_{|k| \leq c_0 K} A_k \cap B_K \right| \leq c_0 K C K^{d-1},
\]
which will be smaller than $K^d \leq |B_K|$ for $c_0$ chosen sufficiently small. Thus, $B_K$ must contain at least one $\vec{v}$ for which $\varphi_k(\vec{v}) \notin \Q$ for all $|k| \leq c_0 K$.
\end{proof}
We now finish the proof of Theorem \ref{main1}.  We let $K$ be a big and fix $v$ as in the claim.  For $|k| \leq c_0 K$, $P_{\mu}(\psi^{\Al}_v(\lambda^k)) \neq 0$ by the consideration at the start of the proof and Lemma \ref{lem::irrational}.

For $k \geq c_0 K$, firstly, $\psi^{\Alf}_v(\lambda^k)=1$ because $v \in \Z[\lambda]$ and $|\lambda^k|_\nu \leq 1$ for all $\nu \in \Sf$, so only the Archimedean contribution matters. We estimate
\[
|\Sigma_<(v\lambda^k)| \leq \sum_{i=0}^{d-1} p |v_i| \sum_{\sigma \in S_<} |\sigma(\lambda)|^{c_0 K} \leq C K c^K
\]
for some $c<1$ and $C<\infty$ depending only on $\lambda$. For $K$ sufficiently big, this will be too close to $0$ to be in $E+\Z$.  This shows that $P_\mu(\psi^{\Al}_v(\lambda^k)) \neq 0$ for $k \geq c_0 K$.

For $k \leq -c_0 K$, we recall that by Proposition \ref{prop::adelefact}, 
\[
\psi^{\Al}_v(\lambda^k)=\psi^{\Ab}_v(-\lambda^k)=e(-\Sigma_>(v\lambda^k))
\]
where $\Sigma_>$, defined analogously to $\Sigma_<$ with $S_>^\infty=\{\sigma: |\sigma(\lambda)|>1\}$ in place of $S_<^\infty$, is the argument of $\psi^{\Ab^\infty}$.  The second equality holds because $\psi^{\Ab^{\mathrm{fin}}}_v(\lambda^k)=1$ for $k\leq-d$. We apply the same reasoning as above to $\Sigma_>(v\lambda^k)$ to see that there can be no $k \leq - c_0 K$ leading to a vanishing coefficient.  We have thus shown that $P_\mu(\psi^{\Al}_v(\lambda^k)) \neq 0$ for all $k \in \Z$, which is in contradiction to Theorem \ref{thm::vanishcoeff}. Thus, $\mu$ and $\lambda$ can not be of maximal entropy.
\end{proof}

\section{Entropy at the non-Archimedean places and proof of Theorem \ref{main2} (without power Fourier decay)} \label{sec::computation}
In this section, we show that we can approximate the Garsia entropy at the finite places from below by the entropy of the finite random variables $Y_n$, which are smoothened versions of $\mu_\beta$.  We then use this construction to prove Theorem \ref{main2},  except for the power Fourier decay on $\Al$. 

As outlined in Subsection \ref{subsec::motivation}, this construction can also be used to calculate the entropy $h_\lambda(\mu)$ for any pair of $\mu$ and $\lambda$ which has no conjugates on the unit circle.  For this,  one has to combine the super-additive sequence $H(Y_n)$ for the non-Archimedean places with one of the methods of \cite{4author} or \cite{bv} for the Archimedean places. 

We let $X_n^\prime:=\beta^{n-1} X_n$ and $M:=M_\beta$.  The main goal of this section is to show that $h_\beta(\mu)=\log M$, which will imply everything else quite quickly.   We define the random variables 
\[
Y_n:=X_n^\prime \mo \beta^n \Z[\beta]=\sum_{j=0}^{n-1} \xi_j \beta^j \mo \beta^n \Z[\beta]
\]
on $G_n$. This has a genuinely different distribution from $X_n^\prime$ because for example $M \beta^{n-1} \in \beta^n \Z[\beta]$ using the minimal equation. One can also define $Y_n$ on $\Ab/\overline{\beta^n\Z[\beta]}$ as
\[
Y_n=X_\beta+\beta^n\overline{\Z[\beta]},
\]
where the bar refers to the closure with respect to the topology of $\Ab$. Because $G_n \cong \overline{\Z[\beta]}/\overline{\beta^n\Z[\beta]}$ canonically, these perspectives are equivalent. So as we said in the introduction,  $Y_n$ is really a version of $X_\beta$ smoothened by adding a small box of proportion $\beta^n$. 

Denote the Shannon entropy of a finitely supported random variable by $H(\cdot)$. The following proposition relates the Shannon entropy of $Y_n$ to the Garsia entropy.
\begin{prop} \label{prop::superadd}
The sequence $H(Y_n)$ is super-additive,  that is,  $H(Y_{n+m}) \geq H(Y_n)+H(Y_m)$ for all $n$ and $m$.  Moreover, 
\[
\lim_{n \to \infty} \frac{H(Y_n)}{n}=h_\beta(\mu).
\]
\end{prop}
Before proving Proposition \ref{prop::superadd}, we state the main result needed in the proof of Theorem \ref{main2}.  The third condition in the proposition below is the same as complete vanishing at level $m$.  

Define the natural projection $\pi: G_n \to G_{n-1}$ with $Y_{n-1}=\pi(Y_n)$ and define for given $x \in G_{n-1}$ the random variable  $Y_m|Y_{m-1}=x$ for $z \in G_n$ by
\[
\PP(Y_m=z|Y_{m-1}=x):=\frac{\PP(Y_m=z, Y_{m-1}=x)}{\PP(Y_{m-1}=x)}.
\]
Recall that the conditional entropy of two random variables is given by $H(X|Y):=\sum_y \PP(Y=y) H(X|Y=y)$ and that the formula $H((X, Y))=H(X)+H(X|Y)$ holds.

\begin{prop} \label{prop::charequi}
The following are equivalent:
\begin{enumerate}
\item $H(Y_n|Y_{n-1})=\log M$
\item For all $x \in G_{n-1}$, $Y_n|Y_{n-1}=x$ is equidistributed on $\pi^{-1}(x)$.
\item For all $\psi \in \widehat{G_n}$ with $\psi \notin \widehat{G_{n-1}}$, 
\[
\EE[\psi(Y_n)]=\prod_{k=0}^{n-1} P_\mu(\psi(\beta^k))=0.
\]
\end{enumerate}
\end{prop}
We now prove Proposition \ref{prop::superadd}.  In the proof, we need the following lemma.
\begin{lemma} \label{lem::infogain}
Let $U, V$ be finitely supported random variables in an additive group and $Z$ be a finitely supported random variable that is a function of $U+V$. Assume that $U$ is independent of both $V$ and $Z$. Then
\[
H(U+V)\geq H(U)+H(Z).
\]
\end{lemma}
\begin{proof}[Proof of Lemma \ref{lem::infogain}]
Note that as $Z$ is determined by $U+V$, the distribution of $U+V$ and $(U+V, Z)$ coincide, which implies
\[
H(U+V)=H(U+V, Z)=H(Z)+H(U+V|Z).
\]
Moreover,  by the independence assumptions,
\[
H(U)=H(U|Z, V)=H(U+V|Z, V) \leq H(U+V|Z).
\]
The last inequality follows from the standard fact that $H(X|Y) \leq H(X)$ for any random variables, which is a consequence of Jensen's formula.
\end{proof}
\begin{proof}[Proof of Proposition \ref{prop::superadd}]
The super-additivity is a consequence of Lemma \ref{lem::infogain} applied to
\begin{align*}
U&:=\sum_{j=n}^{n+m-1} \xi_j \beta^j\mo \beta^{n+m}\Z[\beta],\\ 
V&:=\sum_{j=0}^{n-1} \xi_j \beta^j\mo \beta^{n+m}\Z[\beta],\\
 Z&:=Y_n=U+V \mo \beta^n \Z[\beta].
\end{align*}
Indeed, $U$ has the same distribution as $Y_m$, so $H(U)=H(Y_m)$.

For the second part, recall that $\lim_{n \to \infty} \frac{H(X^\prime_n)}{n}=h_\beta(\mu)$.  We want to show this with $X^\prime_n$ replaced by $Y_n$. 
The inequality $\geq$ is obvious because $H(X^\prime_n) \geq H(Y_n)$ as $Y_n$ is a function of $X^\prime_n$. 

For the converse inequality, we let $k$ be so big that $|\beta|_\nu^k<\eta$ for all $\nu \in S_>$, where $\eta$ is as in Lemma \ref{lem::separation} (applied with the roles of $\beta$ and $\lambda$ switched), that is, such that if $x, y \in \mathrm{supp}(X_n^\prime)$ are such that
\[
|x-y|_{\nu} \leq 2 \eta |\beta|_\nu^n
\]
for all $\nu \in S_>$, then $x=y$. For this choice of $k$, the distribution of $X_n^\prime$ coincides with the one of
\[
\widetilde{Y_n}:=X_n^\prime \mo \beta^{n+k} \Z[\beta].
\]
Adding the independent random variable
\[
\sum_{j=n}^{n+k-1} \xi_j \beta^j \mo \beta^{n+k} \Z[\beta]
\]
can only increase the entropy of $\widetilde{Y_n}$, which implies $H(Y_{n+k}) \geq H(\widetilde{Y}_n)=H(X_n^\prime)$. Therefore,
\[
\lim_{n \to \infty} \frac{H(Y_{n+k})}{n+k} \geq \lim_{n \to \infty} \frac{H(X_n^\prime)}{n+k}=h_\beta(\mu).
\]
\end{proof}
We now show Proposition \ref{prop::charequi}, which is the main step to show maximal entropy in Theorem \ref{main2}.
\begin{proof}[Proof of Proposition \ref{prop::charequi}]
Recall that by definition 
\[
H(Y_n|Y_{n-1})=\sum_{x \in G_{n-1}} \PP(Y_{n-1}=x) H(Y_n|Y_{n-1}=x).
\] 
It is easy to see that for a finitely supported random $U$ variable taking $M$ many values,  $H(U) \leq \log M$, with equality attained if and only if $U$ is equidistributed.  This implies the equivalence of $1.$ and $2.$. 

To see that $2.$ is equivalent to $3.$, we calculate
\[
\EE[\psi(Y_n)]=\EE[\EE(\psi(Y_n)|Y_{n-1})]=\sum_{x \in G_{n-1}} \PP(Y_{n-1}=x) \EE[\psi(Y_n)|Y_{n-1}=x].
\]
We let $\psi_0, \psi_1,  \dots,  \psi_{M-1}$ be a representative system of $\widehat{G_{n}}/\widehat{G_{n-1}}$ with $\psi_0$ the trivial character. The condition in $3.$ is then equivalent to saying that for all $\varphi \in \widehat{G_{n-1}} \subset \widehat{G_{n}}$ and all $j \neq 0$, $\EE[\varphi(Y_n) \psi_j(Y_n)]=0$. Note that uniformly $\varphi(Y_n)=\varphi(x)$ on the set $Y_n|Y_{n-1}=x$, which gives
\[
\EE[\varphi(Y_n) \psi_j(Y_n)]=\sum_{x \in G_{n-1}} \varphi(x) \PP(Y_{n-1}=x) \EE[\psi_j(Y_n)|Y_{n-1}=x]=:\sum_{x \in G_{n-1}} \varphi(x) c_j(x).
\]

Fix an $x$. We now show that for any random variable $U$ on $G_n$ taking values in $\pi^{-1}(x)$, $U$ is equidistributed on $\pi^{-1}(x)$ if and only if $\EE[\psi_j(U)]=0$ for all $j \neq 0$. This is easily seen: We can assume $x=0$, else replace $U$ by $U-y$ for some $y \in \pi^{-1}(x)$. For $x=0$, $U$ takes values in the subgroup $\mathrm{ker}(\pi)$ and is equidistributed if and only if $\EE[\psi(U)]=0$ for all non-trivial $\psi \in \widehat{\mathrm{ker}(\pi)}$ by basic character theory. The collection of $\psi_j$ with $j \neq 0$ restricted to $\mathrm{ker}(\pi)$ is exactly the collection of non-trivial characters in $\widehat{\mathrm{ker}(\pi)}$, showing this fact. We deduce that $Y_n|Y_{n-1}=x$ is equidistributed on $\pi^{-1}(x)$ if and only if $c_j(x)=0$ for all $j=1, \dots, M-1$. 
This implies the direction $2. \Rightarrow 3.$

For the converse direction $3. \Rightarrow 2.$, we note that by basic character theory on $G_{n-1}$, for any fixed $j$,
\[
\forall x \in G_{n-1}: \;  c_j(x)=0 \quad \Leftrightarrow \quad \forall \varphi \in \widehat{G_{n-1}}: \sum_{x \in G_{n-1}} \varphi(x) c_j(x)=0.
\]
Our condition in $3.$ thus implies that $c_j(x)=0$ for any $j$ and any $x$ which per the above implies the claim.
\end{proof}
We are now in the position to prove Theorem \ref{main2}.
\begin{proof}[Proof of Theorem \ref{main2} (without power Fourier decay)]
Recall that $m$ is the integer at which complete vanishing is attained.  Note that by definition, complete vanishing then also holds for all $j \geq m$.  By Proposition \ref{prop::charequi},  this implies that $H(Y_j|Y_{j-1})=\log M$ for all $j \geq m$. Thus, for any $n>m$,
\[
H(Y_n)=H(Y_{m-1})+\sum_{j=m}^{n} H(Y_j|Y_{j-1})=H(Y_{m-1})+(n-m+1) \log M.
\]
Letting $n \to \infty$ and using Proposition \ref{prop::superadd} shows that $h_\beta(\mu)=\log M$. 

Proposition \ref{prop::ifmaxthencont} implies that $\mu_\lambda$ and $\mu_\beta$ are absolutely continuous.  

To show that there is a $C=C(\beta, m)$ such that $\widehat{\mu_\beta}(w)$ for all $|w|_{\Ab}\geq C$,  we note that 
\[
\{v \in \Ab: \psi_v^{\Ab}(\beta^n \Z[\beta])=1\} \rightarrow \widehat{G_n} \quad \mathrm{via} \quad v \mapsto (x  \mapsto \psi_v^{\Ab} (x+\beta^n \Z[\beta])).
\]
For every $w \in \Ab$ we can find an $n$ such that $\psi_w^{\Ab}(\beta^n \Z[\beta])=1$ because characters are continuous.  Complete vanishing at level $m$ also implies complete vanishing at level $n$ for all $n>m$. This means that for such $w$,  $\widehat{\mu_\beta}(w)=0$ unless $w$ is mapped to an element in $\widehat{G_{m-1}}$ under the map above. This is in turn equivalent to saying that $w$ is such that $\psi_w^{\Ab}(\beta^{m-1} \Z[\beta])=1$.  As the set $\{w \in \Ab: \psi_w^{\Ab}(\beta^m \Z[\beta])=1\}$ is compact in $\Ab$,  we can find a $C$ such that no $w$ with $|w|_{\Ab} \geq C$ is contained in it.  
\end{proof}

\section{Proof of Theorem \ref{easycase} and example} \label{sec::easycase}
In this section,  we prove Theorem \ref{easycase} and fill in the details of Example \ref{ex::beta}.  The reason that the bounds are so much better in the case that $(a_1, M)=1$,  where $a_d x^d+\dots+a_1x+M$ is the minimal polynomial of $\beta$,  is that the groups $G_n=\Z[\beta]/\beta^n\Z[\beta]$ are cyclic in this case (easy to see and proven for example as part of Lemma \ref{lem::propGn}). The map $\Phi: G_n \to \Z/M^n\Z$ with $\Phi(1)=1$ is then an isomorphism, which can be seen for example by comparing the order of $1 \in \Z[\beta]$ with the group order.  
\begin{proof}[Proof of Theorem \ref{easycase}]
We know from Theorem \ref{main3} that there is an $n$ such that complete vanishing occurs. Showing this directly would also be quite easy in the setting of Theorem \ref{easycase},  in which the groups $G_n$ are cyclic. This direct proof would not need any of the tools developed in the later sections but only the relatively easy ideas from Claim \ref{claim::limittovanishing} and Claim \ref{enum::newvecs}.  

Fix some big $n$ such that there is complete vanishing at level $n$ and no complete vanishing at level $n-1$.  Our aim is to bound the size of this $n$ in terms of $|E|$,  where
\[
E:=\{\theta \in (0,1): \widehat{\mu}(\theta)=0 \mathrm{\; and \;}\exists n: \theta \in \frac{1}{M^n}\Z\}
\]
are as before the relevant zeroes of $\widehat{\mu}$.  We recall that $P_\mu$ is defined by $P_\mu(e(\theta))=\widehat{\mu}(\theta)$, where $\widehat{\mu}$ is the Fourier transform of $\mu$ on $\R/\Z$.  As the set $e(E)$ consists of roots of $P_\mu$ and $P_\mu$ is a polynomial of degree $|I|$, we have that $|E| \leq |I|$,  where $I$ is the interval $\mu$ is supported on.  

We write $\Phi$ for the isomorphism given by
\[
\Phi: \; G_n \to \Z/M^n\Z,  \quad 1 \in \Z[\beta] \mapsto 1 \in \Z/M^n\Z.
\]
Under this isomorphism, $\beta^j$ is mapped to some $k_j M^j$ with $k_j$ coprime to $M$,  because $\beta^j$ is of order $M^{n-j}$ in $G_n$ and $\Phi$ preserves orders. There is an associated isomorphism $\Psi: \widehat{G_n} \cong \Z/M^n\Z$ with
\[
\psi(y)=e\left(\frac{\Psi(\psi) \Phi(y)}{M^n}\right)
\]
for all $\psi \in \widehat{G_n}, y \in G_n$.  We use here that $\widehat{\Z/M^n\Z}=\Z/M^n\Z$.

Think of $\widehat{G_n}$ as an additive group and let
\[
\langle v, y \rangle_<:=\frac{\Psi(v) \Phi(y)}{M^n}
\]
for $v \in \widehat{G_n}, y \in G_n$ be the dual pairing. Define for $ x \in \bigcup_j \frac{1}{M^j} \Z$ the \textit{weight} of $x$ to be 
\[
w(x):=\min\{n \in \Z_{\geq 0}: M^n x \in \Z\}.
\] 
By the explicit expression for $\beta$ above, we then have that $w(\langle v, \beta^{j+1} \rangle_<)=w(\langle v, \beta^{j} \rangle_<)-1$ for all $v$. Moreover, 
\[
\widehat{G_n} \backslash \widehat{G_{n-1}} = \left\{v \in \widehat{G_n}: w(\langle v, 1 \rangle_<)=n \right\}.
\]
We say that $v \in \widehat{G_n} \backslash \widehat{G_{n-1}}$ \textit{is killed by} $e \in E$ at level $j$ if $\langle v, \beta^j \rangle_< \equiv e \mo \Z$.
\begin{claim} \label{claim::easyhelp}
Let for $e \in E$, 
\[
W(e):=\left\{ v \in \widehat{G_n} \backslash \widehat{G_{n-1}} \; \Big| \; \exists j: \; \langle v, \beta^j \rangle_< =e \mo \Z \right\}
\]
denote the elements killed by $e$. Then $|W(e)|=M^{n-w(e)}$ for any $e \in E$ and for any $e, f \in E$ with $w(e) \leq w(f)$, either $W(f) \subset W(e)$ or $W(f) \cap W(e)=\emptyset$. 
\end{claim}
\begin{proof}
By considering the respective weights, we have that for any $e \in E$, the level $j$ the elements $v \in W(e)$ are killed at must be $j=n-w(e)$.  As for any $v, w$,
\[
\langle v, \beta^j \rangle_< \equiv \langle w, \beta^j \rangle_< \mo \Z \quad  \Leftrightarrow \quad v-w \in \widehat{G_j},
\]
one sees that $W(e)=v+\widehat{G_j}$ for any $v \in W(e)$. As $v+\widehat{G_j} \subset w+\widehat{G_i}$ or $(v+\widehat{G_j}) \cap (w+\widehat{G_i})=\emptyset$ for $j \leq i$ and any $v, w \in \widehat{G_n}$, we deduce the claim. 
\end{proof}
We proceed with the proof of Theorem \ref{easycase}, assuming for now that $M$ is odd.
We want to show that 
\[
|E| \geq (2n-1)(M-1).
\]
By the claim, the most efficient thing to do in order to minimize $|E|$ is to kill as many elements in $\widehat{G_n} \backslash \widehat{G_{n-1}}$ with $e$ of the lowest possible weight, as writing 
\[
W(e)=W(f_1) \cup \dots W(f_r)
\]
with $w(f_i)>w(e)$ will require at least $M$ many elements instead of only one. 

The most efficient way to minimize $E$ is then to take $M-3$ many points of weight one, i.e. of the form $\frac{k}{M}$ with $0<k<M$. To take all $M-1$ is impossible because one then would have complete vanishing at level one,  which implies $n=1$.  That one can not take $M-2$ many points is because 
as $\mu$ is a probability measure, $P_\mu$ is invariant under complex conjugation, implying $E=-E$.  The point $1/2$ is not in $E$ because of the assumption that $M$ is odd.

If some $f, -f$ of weight one are left out, one writes
\[
W(f)=\bigcup_{e \in \frac{1}{M^2} \Z, w(e)=2, W(e) \subset W(f)} W(e),
\]
which is a disjoint union of $M$ many sets. The best one can do here without having total vanishing at level $2$ is to include $2M-2$ points of weight $2$ out of the $2M$ many left. The two $W(-e), W(e)$ left out are then again each a disjoint union of $M$ many $W(c)$ of weight $3$, so by the same reasoning, best is to include $2M-2$ elements $c$ with $W(c) \leq W(e)$. Iterating, we find that the most efficient thing to do is to include $M-3$ many points of weight $1$, take $2M-2$ many of weight $j$ for all $j<n$ and finally, when complete vanishing is to be obtained, $2M$ many of weight $n$. As this is the construction minimizing $|E|$, for any $E$,
\[
|E| \geq M-3 +(n-2)(2M-2)+2M=(2n-1)(M-1).
\]
The bounds are worse in the case that $M$ is even because one can leave out $\frac{1}{2}$ at level $1$ without being forced to also leave out a second element of $E$.  One gets
\[
|E| \geq M-2 + M-2+(n-3)(2M-2)+2M=(2n-2)(M-1).
\]
\end{proof}
To end this section, we fill in the details in Example \ref{ex::beta}.  We recall that in this example,  $\beta$ is a root of 
\[
f_\beta(x)=3x^3+4x^2+3x+5
\]
and that given a $0 \leq k<20$, we wish to classify all measures $\nu$ with maximal entropy supported on $\{0, \dots, k\}$.  The classification will be given via $E_\nu:=\{\theta \in (0,1): \widehat{\nu}(\theta)=0\}$.  We recall $E_1:=\{\frac{j}{5},  1 \leq j \leq 4\}$ and 
\[
E_2^{(1)}:=\left\{\pm \frac{1}{5} \right\} \cup \left\{\pm \frac{1+5j}{25}, 0 \leq j \leq 4 \right\}, \quad E_2^{(2)}:=\left\{\pm \frac{2}{5} \right\} \cup \left\{\pm \frac{2+5j}{25}, 0 \leq j \leq 4 \right\}.
\] 

For $0 \leq k <4$, there are no measures $\nu$ with maximal entropy, which can be seen either by the theorem or because the entropy of $\nu$ is too small. 

For $4 \leq k<12$, complete vanishing must occur on level $m=1$ by Theorem \ref{easycase}, which is equivalent to $E_\nu \supset E_1$. These $\nu$ are exactly those which are equidistributed modulo $5$; the smallest example is the uniform measure on $\{0, 1, 2, 3, 4\}$. 

For $12 \leq k <20$, which means $m \leq 2$ by the theorem, there are of course the measures $\nu$ with $m=1$.  However, there are also other possibilities. To calculate the new options for the zero sets $E$, we consider the group $ G_2=\Z[\beta]/\beta^2 \Z[\beta]$ with the isomorphism
 \[
 \Phi: G_2 \cong \Z /25\Z, \quad 1 \in \Z[\beta] \mapsto 1 \in \Z /25\Z.
 \]
This also induces an isomorphism $\Psi: \widehat{G_2} \cong \Z /25\Z$. In $G_2$,  $5 \equiv -3 \beta \mo \beta^2\Z[\beta]$, so as $5\beta \equiv 25 \equiv 0$ in $G_2$,
 \[
 \beta \equiv 15 \mo \beta^2\Z[\beta].
 \]
 
As before, for $\psi \in \widehat{G_2}$ and $y \in G_2$,
\[
\psi(y)=e\left(\frac{\Phi(y) \Psi(\psi)}{25}\right).
\]
Writing $x=x(\psi)=\Psi(\psi)$,  we then have $\psi(\beta)=e\left(\frac{15x}{25}\right)$.
Thus, if $\pm \frac{1}{5} \in E_\nu$, then all $\psi \in \widehat{G_2}$ with $\frac{15 x}{25}=\pm \frac{1}{5}$
are killed at level one.  To achieve full vanishing at level $2$,  all $\psi \in \widehat{G_2} \backslash \widehat{G_1}$ must be killed at level one or two.  Under the isomorphism, $\psi \in \widehat{G_2} \backslash \widehat{G_1}$ correspond to the $x \in \Z/25\Z$ which are not divisible by $5$.

Assuming complete vanishing at level $2$ and $\pm\frac{1}{5} \in E_\nu$,  one must thus have $\frac{x}{25} \in E_\nu$ for all $x$ such that $\frac{15 x}{25} \neq \pm \frac{1}{5} \mo \Z$. The first genuinely new examples of measures with complete vanishing at level $2$ are thus measures with $E_\nu \supset E_2^{(1)}$. 
Similarly, if one starts by picking $\pm \frac{2}{5}\in E_\nu$ on the first level, one then gets as the other genuinely new possibility for the zero set $E_\nu \supset E_2^{(2)}$.  Picking no element at level $1$ is not an option because $\nu$ will then be a rescaled version of a measure with complete vanishing at level $1$ (and also have support on an interval of length at least $20$ in this particular example).

When one considers $k \geq 20$, one will also have measures appearing for which complete vanishing occurs only at level $3$ or higher. They could in principle all be classified in the manner above.

The example with smallest support of a measure with maximal entropy that is not equidistributed modulo $5$ is the measure $\nu$ supported on $\{0, \dots, 12\}$ with  
 \[
P_\nu(z)=\frac{1}{5 \varphi^2}(z^2+\varphi z+1)(z^{10}+\varphi z^5+1),
 \]
 where $\varphi=-\left[e\left(\frac{2}{5}\right)+e\left(-\frac{2}{5}\right)\right]$ is the golden ratio.  One can check easily that $ E_\nu=E_2^{(2)}$ in this case.  An explicit expression of this measure is given in Example \ref{ex::beta}.

\section{Tools for the case when $\beta$ is not an algebraic integer} \label{sec::notalgint}
Recall that $\beta=\lambda^{-1}$.  In the case when $\beta$ is an algebraic integer, which is equivalent to $\Al=\R^d$, the star roles in the proof of Theorem \ref{main3} and power Fourier decay in Theorem \ref{main2} are played by the lattice $\Z[\beta] \subset \R^d$ and its polar lattice $\Z[\beta]^*$.  These lattices have two properties that are important in the proof.

The first is that to show $\psi_v^{\Al}(\beta^j )=1$ for all $j \geq 0$, it suffices to show $v \in \Z[\beta]^*$, which is implied by $\langle v,  \beta^j \rangle_< \in \Z$ for $0 \leq j<d$.

The second is that the lattice $\Z[\beta]^*$ captures the characters of $G_n$; that is,  $\widehat{G_n}=\lambda^n\Z[\beta]^*/\Z[\beta]^*$.  This follows because $(\beta^n \Z[\beta])^*=\lambda^n \Z[\beta]^*$ and $G_n=\Z[\beta]/\beta^n \Z[\beta]$.

These relatively easy properties are all that is needed to show Theorem \ref{main3} and power decay in Theorem \ref{main2} in the case that $\beta$ is an algebraic integer.  In this section,  we develop analogous tools for the general case.  For this,  we consider 
\[
L:=\sum_{0 \leq j \leq m} \beta^j \Z \subset \R^d,
\]
which is a lattice because it is discrete by the product formula. We show below that there is an $m=m(\beta) \geq d-1$ such that the lattice $L$
and its polar lattice $L^*$ have analogous properties to the two properties of $\Z[\beta]^*$ above.  Proposition \ref{prop::latticefact} below is analogous to the first property (the assumption in the proposition is equivalent to $w \in L^*$) and Lemma \ref{lem::Lstructure} is analogous to the second. 

Recall that $D_\infty$ is the diagonal embedding $D_\infty \colon K \to \R^d$, which we frequently drop, writing for example $\beta \in \R^d$ in place of $D_\infty(\beta)$.  Set $S_==\{\nu: \; |\lambda|_\nu=1\}$ and recall that
\[
\Sf=\{\nu \mathrm{\; non-Archimedean}: \; |\lambda|_\nu<1\}.
\]
\begin{prop} \label{prop::latticefact}
Let $\beta$ be an algebraic number such that $|\sigma(\beta)|>1$ for all Galois conjugates. There is an $m=m(\beta) \geq d-1$ such that if $w \in \R^d$ is such that $\langle w, \beta^j \rangle_< \in \Z$ for every $0 \leq j \leq m$, then $w \in D_\infty(K)$ and 
\[
\forall s \in \Z: \; \psi^{\Al}_w(\beta^s)=\psi^{\Ab}_w(\beta^s)^{-1},
\]
where $w$ is implicitly identified with its preimage $D_\infty^{-1}(w) \in K$.  

Moreover, $\psi_w^{\Alf}(\lambda^j)=1$ for all $j \geq -m$ and $\psi^{\Al}_w( \beta^j)=1$ for all $j \geq 0$.
\end{prop}

For the analogue of the second property,  recall that $G_n \cong \overline{\Z[\beta]}/\beta^n \overline{\Z[\beta]} \subset \Ab/\beta^n \overline{\Z[\beta]}$, where the closure is taken with respect to the topology of $\Ab$.  In particular, this implies that any character of $\Ab$ which is constant on $\beta^n \overline{\Z[\beta]}$ induces a character of $G_n$. Moreover,  by the first point in Proposition \ref{prop::latticefact}, $L^* \subset D_\infty(K)$, so that the mapping in the lemma below makes sense by identifying $v$ with its preimage in $K$.
\begin{lemma} \label{lem::Lstructure}
Let $\beta$ and $L$ be as above.  For any $n \in \N$,
\[
\lambda^n L^*/(\lambda^n L^* \cap L^*) \cong \widehat{G_n},
\]
where the isomorphism is given by $v \mapsto \psi_v^{\Ab}$.  
\end{lemma}
To prove the two results,  we start with two preparatory lemmas.  After this, we introduce some results from algebraic number theory we need in the proof.
\begin{lemma} \label{lem::vinK}
Let $\beta$ be an algebraic number such that $|\sigma(\beta)|>1$ for all Galois conjugates.  If $v \in \R^d$ is such that $\langle v,  \beta^j \rangle_< \in \Q$ for $0 \leq j<d$, then $v \in D_\infty(K)$.
\end{lemma}
\begin{proof} 
We recall from Section \ref{sec::ant} that $\langle x, y \rangle_<=\langle S(x),  y \rangle_{\mathrm{std}}$ for an invertible symmetric matrix $S$,   where $\langle -,- \rangle_{\mathrm{std}}$ denotes the standard scalar product.  We let $A:=(D_\infty(1), \dots, D_\infty(\beta^{d-1}))^T \in \mathrm{Mat}_{d \times d}(\R)$.  By assumption on $v$, we then have that $ASv=(\langle \beta^{j-1}, v \rangle_<)_{1 \leq j \leq d} \in \Q^d$ or,  in other words,  $v \in S^{-1} A^{-1} \Q^d$.

We now note that $A^T$ induces very natural coordinates on $\R^d$: 
\[
D_\infty(K)=\left\{\sum_{0 \leq j<d} q_j D_\infty(\beta^j), q_j \in \Q \right\}=A^T \Q^d.
\]
As 
\[
A S A^T=\left(\langle \beta^{i-1},  \beta^{j-1} \rangle_<\right)_{i, j}=\left(\mathrm{tr}_{K/\Q}(\beta^{i+j-2})\right)_{i, j}
\]
has rational coefficients,  $(A S A^T)^{-1} \Q^d \subset \Q^d$ and thus indeed $S^{-1} A^{-1} \Q^d \subset A^T \Q^d =D_\infty(K)$. 
\end{proof}
Before we prove Proposition \ref{prop::latticefact}, we need some more algebraic number theory.  We let $V_p:=\prod_{\nu|p} K_\nu$, where $\nu|p$ denotes that $\nu$ lies over $p$, or equivalently,  that $K_\nu$ is a field extension of $\Q_p$.  As $V_p \cong K \otimes \Q_p$ (\cite{fouriernf},  Proposition 4-40),  $1,  \beta, \dots, \beta^{d-1}$ is a $\Q_p$-basis of this vector space.  As an additive group,  $V_p$ is self-dual,  with the duality induced by its standard character
\[
\psi^{V_p}(y):=\psi^{\Q_p}(\tr_{V_p}(y))=\prod_{\nu|p} \psi^{K_\nu}(y_\nu),
\]
where
\[
\tr_{V_p}(y):=\sum_{\nu|p} \tr_{K_\nu / \Q_p}(y_\nu).
\]
Furthermore, set $S_{=,p}:=\{\nu|p: |\beta|_\nu=1\}$ and define $V^\prime_p:=\prod_{\nu \in S_{=,p}} K_\nu$,
\[
\tr_{V^\prime_p}(y):=\sum_{\nu \in S_{=,p}} \tr_{K_\nu / \Q_p}(y_\nu)
\]
and $\psi^{V_p^\prime}(y):=\psi^{\Q_p}(\tr_{V_p^\prime}(y))$.  Before we prove Proposition \ref{prop::latticefact}, we need another lemma.  
\begin{lemma} \label{lem::latticehelp}
Let $\beta$ be as above and let $p$ be a prime.   If $w \in K$ is such that $\tr_{V^\prime_p}(w \beta^j) \in \Z_p$ for $0 \leq j<d$, then $\tr_{V^\prime_p}(w \beta^j) \in \Z_p$ for all $j \in \Z$.
\end{lemma}
\begin{proof}
A property we need is the relation of the factorisation of the minimal polynomial $f \in \Z[X]$ of $\lambda$ over $\Q_p$ to the places $\nu|p$.  One has
\[
f(X)=\prod_{\nu|p} f_\nu(X),
\]
where $f_\nu \in \Z_p[X]$ is the minimal polynomial of $\lambda$ embedded into $K_\nu$. One then also has $K_\nu \cong \Q_p[X]/(f_\nu(X))$ (\cite{fouriernf}, Proposition 4-31). One sees directly from the definition that $|\lambda|_\nu>1$ (resp.  $|\lambda|_\nu<1$) if and only if the leading (resp. constant) coefficient of $f_\nu$ lies in $p\Z_p$. 

We set 
\[
f^\prime=\prod_{\nu \in S_{=,p}} f_\nu \in \Z_p[X],
\]
and write $f^\prime(X)=\sum_{0 \leq j\leq s} c_j X^j$ with $c_j \in \Z_p$, where $s\leq d$. By the above, $c_0$ and $c_s$ are invertible in $\Z_p$.  

Denote the diagonal embedding of $y \in K$ into $V_p^{\prime}$ by $D^\prime(y)$.
We note that as all $f_\nu$ with $\nu \in S_{=,p}$ divide $f^\prime$,  
\[
\sum_{0 \leq j \leq s} c_j D^\prime(\beta^j)=0
\]
and hence in particular $\sum_{0 \leq j \leq s} c_j \tr_{V_p^{\prime}}(u \beta^j)=0$ for all $u \in K$.  This shows that for any $n$, 
\[
\tr_{V_p^{\prime}}(w\beta^{n+s})=-c_s^{-1} \sum_{0 \leq j<s} c_j \tr_{V_p^{\prime}}(w\beta^{n+j})
\]
and
\[
\tr_{V_p^{\prime}}(w\beta^{n-1})=-c_0^{-1} \sum_{0 \leq j<s} c_{j+1} \tr_{V_p^{\prime}}(w\beta^{n+j}).
\]
Using these equations iteratively concludes the proof.
\end{proof}
We are now in the position to prove Proposition \ref{prop::latticefact}. We recall that $S_==\{\nu: |\beta|_\nu=1\}$ and that this only contains non-Archimedean places by the assumption on $\beta$.
\begin{proof}[Proof of Proposition \ref{prop::latticefact}]
The first claim, namely that $w \in D_\infty(K)$,  follows from Lemma \ref{lem::vinK}. The main part of the proof is to show  
\begin{equation} \label{eq::factforw}
\forall s \in \Z: \; \psi^{\Al}_w(\beta^s)=\psi^{\Ab}_w(\beta^s)^{-1},
\end{equation}
which we do now.  

For $\beta$ as in Proposition \ref{prop::latticefact},  $S_<$ contains all Archimedean embeddings, so that one can rewrite Proposition \ref{prop::adelefact} as
\begin{equation*} 
e(\langle x, y \rangle_<)=e(\tr_{K/\Q}(xy))=\prod_p \psi^{V_p}(xy)^{-1}
\end{equation*}
for $x,  y \in K$.  From this, or directly using the definition of the trace,  one sees that $\langle x,  y \rangle_< \in \Z$ implies that $\tr_{V_p}(xy) \in \Z_p$ for all $p$.  We deduce that under the assumption on $w$ in Proposition \ref{prop::latticefact},
\begin{equation} \label{eq::wbeta}
\tr_{V_p}(w\beta^j)=\sum_{\nu|p} \tr_{K_\nu / \Q_p}(w\beta^j) \in \Z_p
\end{equation}
for $0 \leq j \leq m$ and all primes $p$. We recall that by Proposition \ref{prop::adelefact},
\begin{equation} \label{eq::charrelation}
\psi^{\Al}_w(\beta^s)\psi^{\Ab}_w(\beta^s)=\prod_{\nu \in S_=} \psi^{K_\nu}_w(-\beta^s)=\prod_p \psi^{\Q_p}(-\tr_{V_p^\prime}(w\beta^s)),
\end{equation}
where the second equality follows directly from the definition. To show the claim, it thus suffices to show that for all $p$ and all $s \in \Z$,  $\tr_{V^\prime_p}(w \beta^s) \in \Z_p$.  

If we knew that $|w \beta^j|_\nu \leq 1$ for all $0 \leq j<d$ and all $\nu$ with $|\beta|_\nu \neq 1$ we could just take $m=d-1$ and would be done, as we argue below Claim \ref{claim::lathelp}.  It will not necessarily be the case that this condition is met.  However, assuming that $m$ as in Proposition \ref{prop::latticefact} is large enough,  the claim below allows us to overcome this.

\begin{claim} \label{claim::lathelp}
If $m=m(\beta)$ is large enough,  then there is a $0 \leq m_0 \leq m-d$ such that for $w^\prime:=\beta^{m_0}$,  $|w^\prime \beta^j|_\nu \leq 1$ for all $0 \leq j<d$ and all $\nu$ with $|\beta|_\nu \neq 1$.
\end{claim}
Before we prove this claim, we show how it implies (\ref{eq::factforw}).  It obviously suffices to show (\ref{eq::factforw}) for $w^\prime$. 

Recalling that for any non-Archimedean $\nu$,  $|x|_\nu \leq 1$ implies that $\tr_{K_\nu/\Q_p}(x) \in \Z_p$, we see from Claim \ref{claim::lathelp} that $\tr_{K_\nu/\Q_p}(w^\prime \beta^j) \in \Z_p$ for all $0 \leq j<d$ and all $\nu$ with $|\beta|_\nu \neq 1$.  By (\ref{eq::wbeta}),  this implies that $\tr_{V^\prime_p}(w^\prime \beta^j) \in \Z_p$ for all primes $p$ and all $0 \leq j <d$.  Lemma \ref{lem::latticehelp} implies that $\tr_{V^\prime_p}(w^\prime \beta^j) \in \Z_p$ then also holds for all $j \in \Z$.   By (\ref{eq::charrelation}), this shows (\ref{eq::factforw}). 
\begin{proof}[Proof of Claim \ref{claim::lathelp}]
We recall from above Lemma  \ref{lem::latticehelp} that $1, \beta,  \dots, \beta^{d-1}$ is a $\Q_p$-basis of $V_p$ and denote its dual basis by $(\beta_j^*)_{0 \leq j <d}$. That is, $\beta_j^*$ are such that 
\[
\tr_{V_p}(\beta_j^* \beta^i)=\delta_{ij}.
\]

We note that the assumption $\tr_{V_p}(w \beta^j) \in \Z_p$ for  $0 \leq j \leq m$ implies that 
\[
\beta^k w \in \sum_{0 \leq j<d} \beta_j^* \Z_p
\]
for all $0 \leq k \leq m-d$. For any $\nu$ with $|\beta|_\nu>1$ (resp. $|\beta|_\nu<1$) there is an $m_\nu \in \N$ such that $|(\beta_j^*)_\nu|_\nu \leq |\beta|_{\nu}^{m_\nu}$ (resp. $|(\beta_j^*)_\nu|_\nu \leq |\beta|_{\nu}^{-m_\nu}$) for all $j <d$. We recall that $S_>=\{\nu: |\lambda|_\nu>1\}$ and set
\[
m:=2d+\max\{ m_\nu: \nu \in S_<^{\mathrm{fin}}\}+\max\{ m_\nu: \nu \in S_> \}.
\]
We let $m_0=\max\{ m_\nu: \nu \in S_> \}$ and let $w^\prime=\beta^{m_0} w$.
The definition of $m_0$ assures that then for any $\nu \in S_>$,
\[
w^\prime\in \beta^{m_0} \sum_{j<d} \beta_j^* \Z_p \subset \OKnu,
\]
where $\OKnu$ is defined by $\OKnu:=\{x \in K_\nu: |x|_\nu \leq 1\}$ for any $\nu$.
Similarly,  for $\nu \in \Sf$,
\[
\beta^d w^\prime=\beta^{-m+m_0+2d} \beta^{m-d}w \in \beta^{-m+m_0+2d} \sum_{j<d} \beta_j^* \Z_p \subset \OKnu.
\]
Thus,  we have indeed $|\beta^j w^\prime|_\nu \leq 1$ for all $\nu$ with $|\beta|_\nu \neq 1$ and for all $0 \leq j<d$.  
\end{proof}
Having shown (\ref{eq::factforw}) and Claim \ref{claim::lathelp},  the final two statements in Proposition \ref{prop::latticefact} follow quite quickly: For any $\nu \in \Sf$, 
\[
|w|_\nu=|\lambda^{m_0}|_\nu |w^\prime|_\nu \leq |\lambda|^{m_0}_\nu.
\]
This implies that $\psi_w^{\Alf}(\lambda^j)=1$ for all $j \geq -m_0$ as $\psi^{K_\nu}(x)=1$ if $|x|_\nu \leq 1 $ for every non-Archimedean $\nu$. For $-m_0 \geq j \geq -m$, this follows from $\psi_w^{\Ab}(\lambda^j)=1$, which we show below,  together with $e(\langle w,  \lambda^j \rangle_<)=1$.

For the final statement,  note that for any $j \geq m_0$ and any $\nu \in S_>$, 
\[
|\beta^j w|_\nu=|\beta|_\nu^{j-m_0} |w^\prime|_\nu \leq 1,
\]
so $\psi^{\Ab}(\beta^j w)=1$ and by (\ref{eq::factforw}) also
\[
\psi^{\Al}(\beta^j w)=(\psi^{\Ab}(\beta^j w))^{-1}=1.
\]

For $0 \leq j<m_0$, 
\[
\psi^{\Al}(\beta^j w)=e(\langle w, \beta^j \rangle_< )=1
\]
because $\psi^{\Alf}(\beta^j w)=1$ for such $j$ as we have shown above. Thus, Proposition \ref{prop::latticefact} is proved.
\end{proof} 
Before we can prove Lemma \ref{lem::Lstructure}, we need one more result exploring properties of the groups $G_n=\Z[\beta]/\beta^n \Z[\beta]$.  We remark that for $f$ the minimal polynomial of $\beta$,  $\Z[\beta] \cong \Z[X]/(f)$.
\begin{lemma} \label{lem::propGn}
Let 
\[
f(X)=\sum_{j=0}^d a_j X^j \in \Z[X]
\]
be irreducible and set
\[
r:=\min\{j: \mathrm{gcd}(a_0,  \dots, a_j)=1\}.
\]
Consider the finite Abelian group $G_n:=\Z[X]/I_n$, where $I_n:=(f, X^n)$ is the ideal generated by $f$ and $X^n$.  This group has $|G_n|=|a_0|^n$ many elements and is generated by $1, X, \dots, X^{r-1}$ as an Abelian group. 
\end{lemma}
\begin{proof}
The group has at most $|a_0|^n$ many elements because 
\[
a_0 X^l=X^l f(X)-(f(X)-a_0)X^l \in (f,  X^{l+1})
\]
for any $l$.  It also has at least that many elements because no polynomial $b_s X^s+b_{s+1} X^{s+1}+\dots$ with $a_0 \nmid b_s$ can be a multiple of $f$.

We show the second property by induction over $n$.  Note that $r \leq d$ because $f$ is irreducible.

Assume  $n=r+1$.  Let $K$ be the subgroup generated by $1, \dots, X^{r-1}$ in $G_n$. We wish to show that $X^r \in K$, for which one only has to show that $a_{j} X^r \in K$ for all $0 \leq j \leq r$ because the $a_j$ have greatest common divisor one.  For $j=r$,  
\[
a_{r} X^r=-a_0 - a_{1} X - \dots - a_{r-1} X^{r-1} \mo I_n.
\]
For other $j$, one can express $a_{j} X^r$ in terms of $a_0 X^{r-j}, \dots,  a_{j-1} X^{r-1}$ in the same way.  Thus indeed $X^r \in K$.

For $n>r+1$,  let $K \subset G_n$ be the subgroup generated by $1, \dots, X^{r-1}$ as above.  By induction,  $X^j \in K+\Z X^{n-1}$ for all $j$.  Thus, for any $0 \leq j<n-1$ there is $c_j \in \Z$ such that $X^j+c_j X^{n-1} \in K$.  

If the numbers $c_j$ were not present,  exactly the same argument as in the case $n=r+1$ would finish the proof. However, the numbers $c_j$ only pose a minor nuisance.  

We have finished the proof if we can show $X^{n-1} \in K$. We note that $\sum_{j=0}^{k} a_{j} X^{n-1-k+j}=0$ in $G_n$ for any $0 \leq k \leq r$. Using this, we see $a_0 X^{n-1} \in K$ and 
\begin{align*}
-a_k X^{n-1}+\sum_{j=0}^{k-1} a_j c_{n-1-k+j} X^{n-1} &= \sum_{j=0}^{k-1} a_{j} X^{n-1-k+j}+\sum_{j=0}^{k-1} a_j c_{n-1-k+j} X^{n-1}\\
&= \sum_{j=0}^{k-1} a_{j} \left(X^{n-1-k+j}+c_{n-1-k+j} X^{n-1}\right) \in K
\end{align*}
for all $0<k \leq r$. 

We let $d_j$ be the smallest common denominator of $a_{0}, \dots, a_{j}$ and set for $0<k \leq r$,
\[
e_{k-1}:=d^{-1}_{k-1}\sum_{j=0}^{k-1}  c_{n-1-k+j} a_{j} \in \Z.
\]
We have just shown that $(a_k-e_{k-1} d_{k-1}) X^{n-1} \in K$ for all $0< k \leq r$.  We note that as
\[
\mathrm{gcd}(a_k-e_{k-1} d_{k-1}, d_{k-1})=\mathrm{gcd}(a_k, d_{k-1})=d_k
\]
for all $k$,  the collection $a_0,  a_1-e_{0} d_{0},  \dots,  a_r-e_{r-1} d_{r-1}$ is still coprime.  This implies that $1$ can be written as a linear combination of them, showing that also $X^{n-1} \in K$.  This concludes the proof.
\end{proof}
 We are now ready to prove Lemma \ref{lem::Lstructure}, which we repeat for convenience.  Recall $L=\sum_{0 \leq j \leq m} \beta^j \Z \subset \R^n$ with $m$ as in Proposition \ref{prop::latticefact}.
\begin{lemma*}
For any $n$, $\lambda^n L^*/(\lambda^n L^* \cap L^*) \cong \widehat{G_n}$, with the isomorphism given by $v \mapsto \psi_v^{\Ab}$.  
\end{lemma*}
\begin{proof}[Proof of Lemma \ref{lem::Lstructure}]
Fix some $n$.  We recall that by Lemma \ref{lem::latticehelp}, $\lambda^n L^* \subset D_\infty(K)$. Implicitly identifying $v \in \lambda^n L^*$ with its preimage in $K$, which in turn embeds into $\Ab$, we can define the map $v \mapsto \psi_v^{\Ab}$. By Proposition \ref{prop::latticefact}, $\psi_v^{\Ab}(\beta^s)=1$ for all $s \geq n$, so $\psi_v^{\Ab}$ is trivial on $\beta^n \Z[\beta] \subset \Ab$. Any character of $\Ab$ that is trivial on $\beta^n \Z[\beta]$ induces a character of $G_n$ in the obvious way. 

We set $\Phi: \lambda^n L^* \to \widehat{G_n}, v \mapsto \psi_v^{\Ab}$. As we have just seen, $\Phi$ is well-defined and it is easy to see that it is a group homomorphism.  We now show that
\[
\mathrm{ker}(\Phi)=\{ v \in \lambda^n L^*: \psi_v^{\Ab}(\beta^s)=1 \; \forall s \geq 0\} = 
\lambda^n L^* \cap L^*.
\]
By Proposition \ref{prop::latticefact}, for  any $v \in L^*$, $\psi_v^{\Ab}(\beta^j)=1$  for all $j \geq 0$, so indeed $\lambda^n L^* \cap L^* \subset  \mathrm{ker}(\Phi)$. For the converse direction,  let $v \in \mathrm{ker}(\Phi)$. We have to show that $\langle v, \beta^j \rangle_< \in \Z$ for $0 \leq j \leq m$.  As $v \in \lambda^n L^*$,  by Proposition \ref{prop::latticefact} $\psi_v^{\Alf}(\beta^j)=1$ for all $j \leq n+m$. Thus,  for any $0 \leq j \leq m$,
\[
e(\langle v,  \beta^j \rangle_<)=\psi_v^{\Al}(\beta^j)=\psi_v^{\Ab}(-\beta^j)=1.
\]
Therefore indeed $\mathrm{ker}(\Phi)=\lambda^n L^* \cap L^*$.\\

It thus only remains to see that $\Phi$ is surjective. For this, it suffices to check that both groups have the same (finite) number of elements. We recall that $|\widehat{G_n}|=|G_n|=M^n$ by Lemma \ref{lem::propGn} applied to the minimal polynomial $f_\beta(X)=\sum_{0 \leq j \leq d} a_j X^d$, $a_0=M$,  of $\beta$. 

To count the elements in the other group, we note that
\[
\left| \lambda^n L^* / \lambda^n L^* \cap L^* \right|=\frac{\mathrm{covol}(\lambda^n L^* \cap L^*)}{\mathrm{covol}(\lambda^n L^*)}.
\]
We set $M_{\infty}:=\prod_{\nu \in \Si} |\lambda|^{-1}_\nu$ and recall that then $|a_d|=\frac{M}{M_\infty}$.
Then, $\mathrm{covol}(\lambda^n L^*)=M_\infty^{-n} \mathrm{covol}(L^*)$ and hence
\[
\mathrm{covol}(\lambda^n L^* \cap L^*)=\left|L^* / \lambda^n L^* \cap L^*\right|\mathrm{covol}(L^*) .
\]
We note that
\[
L^* / \lambda^n L^* \cap L^* \cong (L^* / \lambda^n L^* \cap L^*)^\wedge \cong (\beta^n L +L)/L
\]
because for any lattices $L_1 \subset L_2$, $L_1^*/L_2^* \cong \widehat{L_2/L_1}$ and $(L_1 \cap L_2)^*=L_1^*+L_2^*$. 

To calculate the size of the finite Abelian group $(\beta^n L +L)/L$, we embed it into $\Alf$ and calculate the size of the respective closures with respect to the topology of $\Alf$. By the second statement of Lemma \ref{lem::propGn} applied to the minimal polynomial 
\[
f_\lambda(X)=X^d f_\beta(X^{-1})=a_d+a_{d-1} X+ \dots
\]
 of $\lambda$,  the collection $1, \lambda, \dots, \lambda^{d-1}$ generates $\Z[\lambda]/\lambda^n \Z[\lambda]$ for any $n$. As $\Z \lambda^n \to 0$ as $n \to \infty$ in $\Alf$,  this implies that $\overline{L}=\beta^{m} \overline{\Z[\lambda]}$.  In particular, $\overline{\beta^n L +L}=\beta^n \overline{L}$,
from which we deduce that
\[
(\beta^n L +L)/L \cong \overline{\beta^n L +L}/\overline{L} =\beta^{n+m} \overline{\Z[\lambda]}/\beta^m \overline{\Z[\lambda]} \cong \overline{\Z[\lambda]}/ \lambda^n\overline{\Z[\lambda]} \cong \Z[\lambda]/ \lambda^n \Z[\lambda],
\]
where all closures are still taken in $\Alf$. By the first statement in Lemma \ref{lem::propGn},  $
\left|\Z[\lambda]/ \lambda^n \Z[\lambda]\right|=|a_d|^n$. 

Thus,
\[
\left| \lambda^n L^* / \lambda^n L^* \cap L^* \right|=\frac{\left|L^* / \lambda^n L^* \cap L^*\right|\mathrm{covol}(L^*)}{M_\infty^{-n}\mathrm{covol}(L^*)}=\frac{|a_d|^n}{M_\infty^{-n}}=M^n.
\]
This concludes the proof of Lemma \ref{lem::Lstructure}.
\end{proof}
 
\section{Proof of power Fourier decay in Theorem \ref{main2}} \label{sec::fourierdecay}
\paragraph*{Outline of strategy and proof in the case that $\beta$ is an algebraic integer}
We recall that as always $\beta=\lambda^{-1}$. The key to showing power Fourier decay on $\Al \supset \R^d$ is the proposition below,  of which power decay is a fairly straightforward consequence. The approach of showing power decay via a result like this originates in \cite{feng},  see Proposition 4.7 there.  

\begin{prop} \label{prop::powerdec}
Let $\lambda$ and $\mu$ be as in Theorem \ref{main2}. Assume moreover that the support of $\mu$ is not contained in $q\Z+a$ for any $q \geq 2$ and $a \in \Z$.\\
There is an absolute $m=m(\beta)$ and $\eta>0, \rho <1$ depending only on $\beta$ and $\mu$ such that for any $v \in \Al$ with $|v|_{\Al} \geq 1$ either there is $0 \leq j \leq m$ such that 
\[
|\widehat{\mu_\lambda}(v)| \leq (1-\eta) |\widehat{\mu_\lambda}(\lambda^{j+1} v)|
\]
or there is some $r=r(v) \geq 1$ with
\[
|\widehat{\mu_\lambda}(v)| \leq \rho^r |\widehat{\mu_\lambda}(\lambda^{r} v)|. 
\]
\end{prop}
The value of $m$ that can be taken in the proposition is the one found in Claim \ref{claim::lathelp}.  In the case of $\beta$ being an algebraic integer,  one can take $m=d-1$.  

The main part of this section is devoted to proving this proposition. Before we sketch how it is proved, we quickly show how it implies power Fourier decay.
\begin{proof}[Proof of power Fourier decay in Theorem \ref{main2} assuming Proposition \ref{prop::powerdec}]
The condition that the support of $\mu$ is not contained in any $q\Z+a$ is not very significant as $q X_\lambda+a$ and $X_\lambda$ have the same distribution.  Assume thus that this is satisfied. 

Without loss of generality, $1-\eta \leq \rho^m$, else increase $\rho$. Let $v$ be given and apply the proposition to find $k \geq 1$ with 
\[
|\widehat{\mu_\lambda}(v)| \leq \rho^k |\widehat{\mu_\lambda}(\lambda^k v)|. 
\]
Apply the proposition again to $\lambda^k v$ and iterate this until for some $n$, $|\lambda^n v|_{\Al} \leq 1$. We have thus shown that for any $v$ there is $n \in \N$ such that $|\lambda^n v|_{\Al} \leq 1$ and
\[
|\widehat{\mu_\lambda}(v)| \leq \rho^n |\widehat{\mu_\lambda}(\lambda^n v)| \leq \rho^n \leq |v|_{\Al}^{-\delta n} 
\]
for some $\delta>0$. The last inequality follows because $|\lambda^n v|_{\Al} \leq 1$ implies that $n \geq c \log |v|_{\Al}$ for some absolute $c=c(\lambda)>0$.
\end{proof}

We now explain the strategy of the proof of Proposition \ref{prop::powerdec}.  For this,  we fix a $v=(v_f, v_\infty) \in \Al$ failing the first condition.  We aim to show that it must satisfy the second.
For this,  we fix $\eps>0$ and find $w \in K$ with several properties.  The first is that $w$ approximates $v$ in $\Alf$ in the sense that $\psi_{v_f}^{\Alf}(\lambda^j)=\psi_w^{\Alf}(\lambda^j)$ for all $j \geq 0$ and in $\Ali \cong \R^d$ in the sense that $|v_\infty-w|_{\R^d}<\eps$.  The second is that $w \in \CC$,  where we recall that 
\[
\CC=\left\{x \in K \; \Big| \; \forall s \in \Z: \psi^{\Al}_x(\lambda^s)=\psi^{\Ab}_x(\lambda^s)^{-1} \right\}.
\]
As $w \in \CC$,  Theorem \ref{thm::vanishcoeff} is applicable for $w$, implying that there is $r \in \Z$ such that $P_\mu(\psi_w^{\Al}(\lambda^r))=0$.  We will also show that $r>0$.  

We now show how finding such a $w$ allows us to conclude the proof.

\begin{claim} \label{claim::approxfinish}
Let $\eps=\eps(\lambda, \mu)>0$ be sufficiently small.  If for a given $v \in \Al$ of sufficiently big norm $|v|_{\Al}$ there is a $w$ with the properties above, then the second point of Proposition \ref{prop::powerdec} holds. That is, there is $r>0$ such that $
|\widehat{\mu_\lambda}(v)| \leq \rho^r |\widehat{\mu_\lambda}(\lambda^{r} v)|$.
\end{claim}
\begin{proof}
By Theorem \ref{thm::vanishcoeff}, either $w=0$ (in which case $v$ is of bounded norm,  concluding the proof), or there is $r$ such that $P_\mu(\psi_w^{\Al}(\lambda^r))=0$. Fix this $r$, recall that we assume $r>0$ and note that
\[
\left| \psi^{\Al}_{w}(\lambda^r)-\psi^{\Al}_v(\lambda^r) \right|=\left| \psi^{\Alf}_{v_f}(\lambda^r) e(\langle w-v_\infty, \lambda^r \rangle_<) \right|=O\left(\left|\lambda^r(w-v_\infty)\right|_{\R^d} \right).
\]
Thus,  arguing with Lipschitz continuity of $P_\mu$ (which we make precise in Claim \ref{claim::meas}) and setting $\rho:=\max\{|\lambda u|_{\R^d},  |u|_{\R^d}=1\}<1$,
\[
|P_\mu(\psi^{\Al}_{v}(\lambda^r)|=\left|P_\mu(\psi^{\Al}_{w}(\lambda^r)\right|+O\left(\left|\lambda^r(v_\infty-w)\right|_{\R^d} \right)=O(\rho^r \eps)\leq \rho^{r},
\]
where the last inequality holds if we choose $\eps$ small enough.  This implies the second condition because $\widehat{\mu_\lambda}(v)=\prod_{j \geq 0} P_\mu(\psi_v^{\Al}(\lambda^j))$.
\end{proof}
We now fill in the details of how such a $w$ can be found, which is the main part of the proof.  Before we show how this is done in the general case, we state and prove the much simpler result in the case that $\beta$ is an algebraic integer.  Together with the preceding pages,  this then gives the complete proof of power decay in the algebraic integer case.  

We recall that if $\beta$ is an algebraic integer,  then $\Z[\beta]$ is a lattice in $\R^d$ and recall that $\Z[\beta]^*$ is its polar lattice with respect to the bilinear form $\langle -, - \rangle_<$.
\begin{lemma} \label{lem::geom}
Assume that $\beta$ is an algebraic integer.  For every $\eps>0$ there is $\eta>0$ such that: If $v \in \R^d$ is such that $|e(\langle v, \lambda^j \rangle_<)-1| \leq \eta$ for all $0 \leq j<d$, then there is $w \in \Z[\beta]^*$ such that $|v-w|_{\R^d} \leq \eps$.
\end{lemma}

\begin{proof}
The condition $|e(\theta)-1| \leq \eta$ is the same as saying that $\theta$ is very close to an integer. If we only knew that $\langle v, \lambda \rangle_<$ is close to an integer, we could not say much because the kernel of $\langle \cdot, \lambda  \rangle_<$ is just some hyperplane.  However, if we look at all of the first $d$ coefficients, this changes drastically. This is because the real $d \times d$ matrix $A$ with 
\[
Av=(\langle v,  \lambda ^j \rangle_<)_{0 \leq j<d}
\]
for any $v \in \R^d$ is invertible.  Explicitly, $A:= (D(1), \dots, D(\lambda ^{d-1}))^T \cdot S$, where $D$ is the diagonal embedding of $K$ into $\R^d \cong \Al$ and $S$ is the invertible symmetric matrix relating $\langle \cdot,  \cdot \rangle_<$ to the standard scalar product (see Section \ref{sec::ant}). 

As $A^{-1}$ is continuous,  if for $v \in \R^d$, $|\langle v,  \lambda ^j \rangle_<-k_j|<\eta$ for all $0 \leq j < d$ for some integers $\vec{k}=(k_j) \in \Z^d$, then necessarily $| v - A^{-1} \vec{k} |_{\R^d} < C \eta$
for some $C$ depending only on $\lambda$.  One can then set $w:=A^{-1} \vec{k}$.

Showing $A^{-1}\Z^d \subset \Z[\beta]^*$ concludes the proof. To show this, it is enough to show that for any $u \in \R^d$,  $\langle u, \beta^j \rangle_< \in \Z$ for $-d <j \leq 0$ (which is equivalent to $u \in A^{-1} \Z^d$) implies that $\langle u, \beta^j \rangle_< \in \Z$ for all $j>-d$.  

For this,  we recall that if $\beta$ is an algebraic integer, then it has minimal equation $\beta^d=-\sum_{0 \leq j<d} a_j \beta^{j}$ with $a_j \in \Z$.  This implies that for any $n$, $\beta^n=-\sum_{0 \leq j<d} a_j \beta^{n-d+j}$.  Applying this equation iteratively for $n=1, 2, 3, \dots$ concludes the proof of the lemma. 
\end{proof}
\paragraph*{Proof in the general case}
We now state the necessary results in the general case.  The lemma below is relatively easy to show.
\begin{lemma} \label{lem::nonAplaces}
Let $v=(v_f, v_\infty) \in \Al$ and $m \geq 1$ be given. There is an $u \in \CC$ such that $\psi^{\Al}_u(\beta^j)=1$ for all $j\geq 0$ and $\psi^{\Alf}_{v_f}(\lambda^r) =\psi^{\Alf}_{u}(\lambda^r)$ for all $r \geq -m$. 
\end{lemma}
Showing the result below is the main part of the proof and relies on the material in Section \ref{sec::notalgint}.  Studying the easier proof of Lemma \ref{lem::geom} before the proof of Lemma \ref{lem::powerdecgeom} is recommended as the main idea comes across more clearly in the less technical setting. 

\begin{lemma} \label{lem::powerdecgeom}
Let $m$ be as in Proposition \ref{prop::latticefact}.  For any $\eps>0$ there is $\eta>0$ such that:

Let $v=(v_f, v_\infty) \in \Al$ be such that $|\psi_v^{\Al}(\beta^j)-1| \leq \eta$ for all $0 \leq j \leq m$ and let $u$ be as in Lemma \ref{lem::nonAplaces}.  Then, there is $y \in \CC$ such that $|v_\infty-u-y|_{\R^d}< \eps$.  Moreover,  $\psi_y^{\Al}(\beta^j)=1$ for all $j \geq 0$ and $\psi_y^{\Alf}(\beta^s)=1$ for all $s \leq m$.  
\end{lemma}
Once these two results are shown, one can set $w:=\beta^m(u+y)$ and use Claim \ref{claim::approxfinish} to finish the proof of Proposition \ref{prop::powerdec}. 

The proof of Lemma \ref{lem::nonAplaces} relies on the following corollary of the adelic strong approximation theorem (\cite{strongapprox},  Chapter II.15).  It is a generalisation of the Chinese remainder theorem (for $K=\Q$, it is exactly equivalent to the CRT).
\begin{theorem} \label{thm::strongapprox}
Let $S$ be a finite set of non-Archimedean places of $K$,  let $v_\nu \in K_\nu$ and $\eps_\nu>0$ for $\nu \in S$.  Then there is $u \in K$ such that $|u-v_\nu|_\nu \leq \eps_\nu$ for $\nu \in S$ and $|u|_\nu \leq 1$ for all non-Archimedean $\nu \notin S$ (without condition on the Archimedean places).
\end{theorem}
\begin{proof}[Proof of Lemma \ref{lem::nonAplaces}]
By Theorem \ref{thm::strongapprox}, there is $u \in K$ such that $|u-v_\nu| \leq |\lambda|^m$ for all $\nu \in \Sf$ and $|u|_\nu \leq 1$ for all $\nu \notin S_<$.  Recall that for every non-Archimedean $\nu$,  $|x|_\nu \leq 1$ implies $\psi^{K_\nu}(x)=1$. As $|\lambda^r u|_\nu \leq |\lambda^{-m} u|_\nu \leq 1$ for $\lambda \in \Sf$,  we see that $\psi^{\Alf}_u(\lambda^r)=\psi^{\Alf}_v(\lambda^r)$ for all $r \geq -m$.  Moreover,  as $|\beta^s u|_\nu=|u|_\nu \leq 1$ for all $\nu \in S_=$ and all $s$, we also see that $u \in \CC$ (implicitly using Proposition \ref{prop::adelefact}).  Lastly,  as $|\beta^j u|_\nu \leq |u|_\nu \leq 1$ for all $j \geq 0$ and all $\nu \in S_>$,  $\psi^{\Ab}_u(\beta^j)=1$.  
\end{proof}
We now prove Lemma \ref{lem::powerdecgeom}, which is the last missing piece in the proof of power decay.
\begin{proof}[Proof of Lemma \ref{lem::powerdecgeom}]
Set $x:=v_\infty-u \in \R^d$ and note that by the properties of $u$,  
\[
e(-\langle u,  \beta^j \rangle_<)=\psi^{\Al}_{u}(-\beta^j)\psi^{\Alf}_{u}(\beta^j)=\psi^{\Ab}_{u}(\beta^j) \psi^{\Alf}_{u}(\beta^j)=\psi^{\Alf}_{v}(\beta^j)
\]
and hence
\[
e(\langle x,  \beta^j \rangle_<)=\psi^{\Al}_{v}(\beta^j)
\]
for all $0 \leq j \leq m$.  By the assumption on $\psi^{\Al}_v(\beta^j)$ in the lemma,  therefore $\Vert \langle x,  \beta^j \rangle_< \Vert= O(\eta)$ for $0 \leq j \leq m$,  where $ \Vert \theta  \Vert$ for $\theta \in \R$ is the distance to the closest integer.
Arguing with the invertibility of an appropriate matrix as in the proof of Lemma \ref{lem::geom}, we find an element $y \in \R^d$ with $\langle y,  \beta^j \rangle_< \in \Z$ for $0 \leq j<d$ and 
$| x -y |_{\R^d}=O_{\beta}(\eta)$. This implies $ \Vert \langle y,  \beta^j  \rangle_<  \Vert=O_{\beta, m}(\eta)$ for all $0 \leq j \leq m$.  All of the claimed properties of $y$ in Lemma \ref{lem::powerdecgeom} now follow from Proposition \ref{prop::latticefact} once the claim below is established. 
\end{proof}
\begin{claim}
For small enough $\delta=\delta(\beta, m)>0$,  the condition $  \Vert \langle y,  \beta^j  \rangle_<  \Vert \leq \delta$ for $0 \leq j \leq m$ extends $\langle y,  \beta^j \rangle_< \in \Z$ from the range $0 \leq j<d$ to the range $0 \leq j \leq m$.
\end{claim}
\begin{proof}
By Lemma \ref{lem::latticehelp},  the vector $y$ lies in (the image of the diagonal embedding of) $K$.  Recall from above Lemma \ref{lem::latticehelp} that for any prime $p$, $V_p:=\prod_{\nu|p} K_\nu$,  
\[
\mathrm{tr}_{V_p}(x):=\sum_{\nu|p} \mathrm{tr}_{K_\nu/\Q_p}(x)
\]
and $\psi^{V_p}(x)=\psi^{\Q_p}(\mathrm{tr}_{V_p}(x))$. Also recall from above (\ref{eq::wbeta}) on page \pageref{eq::wbeta} that $\langle y,  \beta^j \rangle_< \in \Z$ implies $\mathrm{tr}_{V_p}(\beta^j w) \in \Z_p$ for every prime $p$ and for all $0 \leq j<d$.  

We now consider two cases of primes $p$ separately: Let $a_d X^d+ \dots+a_0$ denote the minimal equation of $\beta$ over $\Z$.  The first case is the one of primes $p$ who do not divide $a_d$.  As $a_d$ is then invertible in $\Z_p$,  using this equation iteratively as at the end of the proof of Lemma \ref{lem::geom},  we see $\mathrm{tr}_{V_p}(\beta^j y) \in \Z_p$ for all $j \geq 0$.  

For the finitely many primes $p$ dividing $a_d$, we note that in this situation, $y \in U_p:=\sum_{j<d} \beta_j^* \Z_p$, where $(\beta_j^*)$ is the dual basis of $(\beta^j)_{j<d}$ in $V_p$. We argued above Lemma \ref{lem::latticehelp} why this dual basis exists. As $U_p$ is compact, there is an absolute $k=k(\beta, p, m)$ such that $\tr_{V_p}(\beta^j y) \in p^{-k}\Z_p$ for all $0 \leq j \leq m $.  There are only finitely many exceptional primes, so the values of 
\[
e(\langle y, \beta^s \rangle_<)= \prod_{\nu \, \mathrm{fin.}} \psi^{K_\nu}(y\beta^s)^{-1}
\]
with $0 \leq s \leq m$ and $y$ ranging over all vectors with $(\langle y,  \beta^j \rangle_<)_{0\leq j<d} \in \Z^d$ are discrete in the unit circle.  Choosing $\delta$ small in terms of the distance between the distinct possible values of $e(\langle y, \beta^s \rangle_<)$ finishes the proof of the claim and thus also the proof of Lemma \ref{lem::powerdecgeom}.
\end{proof} 

We conclude this section by showing the claim below, which was used in the proof of Claim \ref{claim::approxfinish}.
\begin{claim} \label{claim::meas}
Let $\mu$ be a finitely supported measure on $\Z$ whose support is not contained in any $q\Z+a$ for $q \geq 2$ and $a \in \Z$.  For any $\eta^\prime>0$ there is $\eta>0$ such that $|P_\mu(z)| \geq 1-\eta$ implies $|z-1| \leq \eta^\prime$ for all $z \in \mathbb{S}^1$.
\end{claim}
\begin{proof}
The function $|P_\mu|$ is continuous and $|P_\mu|\leq 1$.  To show the claim, it thus suffices to show that $|P_\mu|$ attains its maximum only at $|P_\mu(1)|=1$. 

This follows from the condition on the support of $\mu$:  If $0< x<1$ were such that 
\[
1=|P_\mu(e(x))|=\left| \sum_{n \in \mathrm{supp}(\mu)} \mu(n) e(xn) \right|,
\]
then this necessarily implies that $xn_1=xn_2 \mo \Z$ for all $n_1, n_2$ in the support of $\mu$.  Hence,  $x=\frac{p}{q}$ must be rational, with $q$ dividing all integers in the support of $\mu$.  Contradiction.
\end{proof}

\section{Proof of Theorem \ref{main3}}\label{sec::main3}
We start with an outline of the proof of Theorem \ref{main3} and show how it follows from Claim \ref{claim::limittovanishing}, Claim \ref{enum::newvecs},  Theorem \ref{cor::polarlattice} and Proposition \ref{prop::boundedn0}. We then prove these four results.

The key object in the proof are lattices in $\R^d$. In the case when $\beta$ is an algebraic integer,  we would consider the lattice $\Z[\beta]^* \subset \R^d$.  In the case in which $\beta$ is not an algebraic integer,  
we consider the lattice $L^* \subset \R^d$, which is the polar lattice of
\[
L:=\sum_{0 \leq j \leq m} \beta^j \Z.
\]
Here,  $m=m(\beta) \geq d-1$ is the integer from Proposition \ref{prop::latticefact}. Our results in Section \ref{sec::notalgint} show that this lattice $L^*$ has properties analogous to the ones $\Z[\beta]^*$ has in the case when $\beta$ is an algebraic integer.  Among other things,  we have shown that $\psi^{\Al}_v(\beta^j)=1$ for all $j \geq 0$ and all $v \in L^*$ and that Theorem \ref{thm::vanishcoeff} applies to elements in $\lambda^n L^*$ for any $n \geq 0$. We then know that for any $v \in L^*$ there is (a smallest) $r=r(v)>0$ with $P_\mu(\psi^{\Al}_v(\lambda^r))=0$.  Ranging over $L^*$, this can of course be arbitrarily big - just consider elements in $\beta^nL^* \cap L^*$ for big $n$.  What needs to be shown in Theorem \ref{main3} is that this $r$ is uniformly bounded for $v \in L^*$, $v \notin \beta L^*$.

We say that $v \in \bigcup_{n \geq 0} \lambda^n L^*$ \textit{has not yet vanished} if
\[
\forall j \geq 0: \quad P_\mu(\psi^{\Al}_v(\beta^j )) \neq 0.
\]
Our aim is to find an absolute $n_0$ such that any $v$ which has not yet vanished is contained in $\bigcup_{n \leq n_0} \lambda^n L^*$.

The first ingredient in the proof is a result showing that if we can find a $v$ which has not yet vanished of small enough norm,  then maximal entropy cannot be achieved.
\begin{claim} \label{claim::limittovanishing} 
There is an $\eps_0>0$ depending only on $\beta$ and the length of the interval $\mu$ is supported on such that if there is a non-zero $v \in \bigcup_{n \geq 0} \lambda^n L^*$ which has not yet vanished and for which $| v |_{\R^d} <\eps_0$, then $h_\lambda(\mu)<\log M_\lambda$.
\end{claim}

The starting point of the proof of Theorem \ref{main3} is to fix a big $n$ and to assume for contradiction that there is a $v \in \lambda^n L^*$ which has not yet vanished.  Whether a point has vanished only depends on the zeroes of $P_\mu$ and among those,  only on zeroes of the form $e(\frac{k}{M^n})$ for $k,  n \in \N$, where $M=M_\lambda$ is the Mahler measure.  There are only finitely many such zeroes,  which we denote (as a subset of $(0,1)$ instead of the unit circle) by
\[
E:=\{\theta \in (0,1): \widehat{\mu}(\theta)=0 \mathrm{\; and \;}\exists n: x \in \frac{1}{M^n}\Z\}.
\]

The claim below now gives a method to generate more $w$ which have not yet vanished from $v$.
\begin{claim}\label{enum::newvecs}
Let $s=s(M, E)$ be the minimal $s$ such that $M^sE \subset \Z$.  Let $v \in \bigcup_{n \geq 0} \lambda^n L^*$ be given which has not yet vanished. Then:
\begin{enumerate} 
\item For $y \in L^*$, $w=v+y$ has not yet vanished.
\item For any $j \geq 0$, $w=\beta^j v$ has not yet vanished.
\item For any $m \in \Z$, $w=(1+M^s m) v$ has not yet vanished.
\end{enumerate}
\end{claim}

One now fixes a $v$ which has not yet vanished and considers the lattices
\[
L_1:=v\Z+L^* \mathrm{\quad and \quad} L_2:=vM^s\Z+L^*
\]
in $\R^d$.  If $L_2$ were $\eps_0$-dense in $L_1$,  then there is a $w=v+mM^sv+y$ with $m \in \Z$,  $y \in L^*$ of norm smaller than $\eps_0$.  The last two claims imply that this is in contradiction to maximal entropy unless $w=0$, which will be a case that is easy to handle.  To get something meaningful out of the case that $L_2$ is not $\eps_0$-dense in $L_1$, we use the following result from the geometry of numbers.
\begin{theorem} \label{cor::polarlattice}
Let $K_1 \subset K_2 \subset \R^d$ be lattices; note that then $K_2^* \subset K_1^*$. There is a constant $C$ only depending on $d$ such that if 
\[
\min\left(| v |_{\R^d}: v \in K_2, v \notin K_1 \right) \geq C \eps^{-1},
\]
then $K_2^*$ is $\eps$-dense in $K_1^*$.
\end{theorem}
Applying this result with $K_1=L_1^*$ and $K_2=L_2^*$ gives an absolute $C_0=C_0(\beta, \mu)$ such that we derive a contradiction to maximal entropy unless $L_2^* \backslash L_1^*$ contains a vector $y$ with norm $|y|_{\R^d} \leq C_0$.  It is easy to see that 
\[
L_2^* \backslash L_1^*=\left\{y \in L: \langle v, y \rangle_< \in \frac{1}{M^s} \Z \mathrm{\; but \;} \langle v, y \rangle_< \notin \Z \right\}.
\]
As $L$ is a lattice,  the set $U:=B_{C_0} \cap L$ is finite.

The final ingredient in the proof is the proposition below.  It states that while for a given $v$ there may well be such an exceptional vector $y \in U$,  this will not be true for all $\beta^m v$ with $0 \leq m  \leq n_0$ for some constant $n_0$ independent of $n$.
\begin{prop} \label{prop::boundedn0}
Let $\beta$ as in Theorem \ref{main3} and $L$ be as above. Let $U \subset L$ be a finite set and let $s \in \N$. Say that $v \in \bigcup_{n \geq 0} \lambda^n L^*$ has property $(\star)$ if there is $y \in U$ such that $\langle v, y \rangle_< \in \frac{1}{M^s}\Z$ but $\langle v, y \rangle_< \notin \Z$. 

There is $n_0=n_0(U, \beta, s)$ such that for any $v \in \bigcup_{n \geq 0} \lambda^n L^*$ there is an $ m \leq n_0$ such that $\beta^m v$ does not have property $(\star)$. 
\end{prop}
We now conclude the proof of Theorem \ref{main3} assuming these results.
\begin{proof}[Proof of Theorem \ref{main3} assuming Claim \ref{claim::limittovanishing}, Claim \ref{enum::newvecs},  Theorem \ref{cor::polarlattice} and Proposition \ref{prop::boundedn0}]
Assume that maximal entropy is attained.  We aim to show that there is complete vanishing at some level $m \leq n_0$.  For this,  it suffices to show by Lemma \ref{lem::Lstructure} that any $v \in \bigcup_{n \geq 0} \lambda^n L^*$ which has not yet vanished must be contained in $\lambda^{m} L^*$ for some $m \leq n_0$.  

Fix a big $n$, let $v \in \lambda^n L^*$ be given that has not yet vanished,  and consider the lattices $L_1$ and $L_2$ defined above.  As we deduced above from Claim \ref{claim::limittovanishing} and Claim \ref{enum::newvecs}, if $L_1$ is $\eps_0$-dense in $L_2$,  then there are $m \in \Z$ and $y \in L^*$ such that $(1+mM^s)v+y=0$.  In this case,  $(1+M^sm)v \in L^*$ maps to the trivial homomorphism in $\widehat{G_n}$ under the isomorphism in Lemma \ref{lem::Lstructure}. As $(1+M^sm)$ is coprime to the group order of $G_n$, this already implies that $v$ maps to the trivial homomorphism, or, equivalently, $v \in L^*$.  Therefore,  if $v \notin L^*$,  then $L_1$ is not $\eps_0$-dense in $L_2$. 

Applying Theorem \ref{cor::polarlattice} as we did above,   and recalling the definition of property $(\star)$ from Proposition \ref{prop::boundedn0},  we see that any $v \in \lambda^n L^*$ which has not yet vanished has property $(\star)$, unless $v \in L^*$.  

We now apply Proposition \ref{prop::boundedn0}. Given $v$, we get an $m \leq n_0$  such that $\beta^m v$ does not have property $(\star)$.  As $\beta^m v$ has also not yet vanished by Claim \ref{enum::newvecs},  $\beta^{m} v \in L^*$ by the above.  Reformulating this as $v \in \lambda^m L^*$ concludes the proof.
\end{proof}
We now prove the results above, starting with the two claims.
\begin{proof}[Proof of Claim \ref{claim::limittovanishing}]
As $\lambda$ acts contractively on $\R^d$,  $|\langle \lambda^j,  v \rangle_<|\leq C |v|_{\R^d} \leq C\eps_0$,  where $C$ is a constant depending only on $\lambda$.  
There is an absolute $\eps_1>0$ depending only on the length of the interval $\mu$ is supported on such that $P_\mu(e(\theta)) \neq 0$ for $|\theta|<\eps_1$. This is because $P_\mu(1)=1$ and its derivative is bounded in terms of its degree as a polynomial, which equals the length of the interval $\mu$ is supported on.  Picking $\eps_0<\eps_1/C$ implies that $P_\mu(e(\langle v, \lambda^j\rangle_<)) \neq 0$ for all $j \geq 0$.  

We recall that $e(\langle v, \lambda^j\rangle_<)=\psi^{\Al}_v(\lambda^j )$ for $j \geq 0$ by Proposition \ref{prop::latticefact},  so this implies that $P_\mu(\psi^{\Al}_v(\lambda^s))) \neq 0$ for all $s \geq 0$.  By the condition that $v$ has not yet vanished, this is also true for $s \leq 0$.  As $v \in \CC$ by Proposition \ref{prop::latticefact},  this is in contradiction to Theorem \ref{thm::vanishcoeff} unless $h_\lambda(\mu)<\log M_\lambda$.
\end{proof}
We now show the second claim.
\begin{proof}[Proof of Claim \ref{enum::newvecs}]
The first statement follows because $\psi^{\Al}_y(\beta^j )=1$ for all $j \geq 0$ by Proposition \ref{prop::latticefact}. The second statement is obvious.  

For the third statement, let $v \in \lambda^n L^*$ be given. For $k \geq n$,  $\psi^{\Al}_v(\beta^k)=1$ implying $\psi^{\Al}_w(\beta^k)=1$, so such $k$ can not lead to vanishing.  For $k<n$,  $\psi^{\Alf}(\beta^k)=1$ by Proposition \ref{prop::latticefact}, so by definition of $E$,  $P_\mu(\psi^{\Al}_w(\beta^k))=0$ if and only if $\langle v, \beta^k \rangle_< \in E+\Z$.

We distinguish between two cases.  If $k<n$ is such that $\langle v, \beta^k \rangle_< \in \frac{1}{M^s} \Z$,  then
\[
\langle w, \beta^k \rangle_<
\equiv \langle v, \beta^k \rangle_< \mo \Z,
\]
which implies $\langle w, \beta^k \rangle_< \notin E+\Z$ because $v$ has not yet vanished.  For $k<n$ such that $\langle v, \beta^k \rangle_< \notin \frac{1}{M^s} \Z$,  also $\langle w, \beta^k \rangle_< \notin \frac{1}{M^s} \Z$, so once again $\langle w, \beta^k \rangle_< \notin E+\Z$. 
\end{proof}
We now state the theorem from the geometry of numbers we need to prove Theorem \ref{cor::polarlattice}. We recall from Section \ref{sec::ant} that the bilinear form $\langle -, - \rangle_<$ is related to the standard scalar product $\langle -, - \rangle_{\mathrm{std}}$ on $\R^d$ via 
\[
\langle x, y \rangle_<=\langle S(x),  y \rangle_{\mathrm{std}},
\]
where $S: \R^d \to \R^d$ is an invertible linear map with $1 \leq |S(x)|_{\R^d} \leq 2$ for all $|x|_{\R^d}=1$. 

We recall that any lattice $K \subset \R^d$ can be written as $K=e_1\Z+\dots+e_d\Z$ for a basis $e_1, \dots, e_d$ of $\R^d$. The dual basis of $e_1, \dots, e_d$ with respect to the bilinear form $\langle -, - \rangle_<$ is the basis $e_1^*, \dots, e_d^*$ such that $\langle e_i^*, e_j \rangle_<=\delta_{ij}$. We define the dual basis $e_{1, \mathrm{std}}^*, \dots, e_{d, \mathrm{std}}^*$ with respect to the standard scalar product in the same way.  We note that as then $S(e_i^*)=e_{i, \mathrm{std}}^*$,  their norms differ at most by a bounded constant.

The $j$-th \textit{successive minimum} of $K$ is defined to be
\[
\lambda_j:=\min\{r: \mathrm{\; there \; are \; independent\;} x_1, \dots, x_j \in K \cap B_r \},
\]
where $B_r=\{x: | x |_{\R^d} \leq r\}$.
\begin{theorem}[Corollary VIII.5.2 in \cite{cassels}] \label{thm::polar}
Let $K \subset \R^d$ be a lattice. There is a basis $e_1, \dots, e_d$ of $K$ such that for all $1 \leq j \leq d$,
\[
\lambda_j \leq | e_j |_{\R^d} \leq j \lambda_j
\]
and 
\[
| e_{j, \mathrm{std}}^* |_{\R^d} \leq C_0 | e_j |_{\R^d}^{-1} ,
\]
where $C_0=C_0(d)$ is a constant depending only on $d$.
\end{theorem}
We now prove Theorem \ref{cor::polarlattice} using this theorem.  We note that we can replace $e_{j, \mathrm{std}}^*$ in the above by $e_j^*$,  changing the constant by a factor of two.
\begin{proof}[Proof of Theorem \ref{cor::polarlattice}]
We let $e_1, \dots, e_d$ be the basis of $K_2$ as in Theorem \ref{thm::polar}.  Setting $C=d C_0$ in Theorem \ref{cor::polarlattice}, we find that for any $j$, either $e_j \in K_1$ or 
\[
| e_j^* |_{\R^d} \leq C_0 | e_j |_{\R^d}^{-1}\leq \frac{C_0}{C \eps^{-1}}=\frac{\eps}{d}. 
\]
Set $S:=\{j: e_j \in K_1\}$.  Let $x \in K_1^*$ be given and write $x=\sum_{j=1}^d a_j e_j^*$
with $a_j \in \R$. Note that for $i \in S$,  $a_i=\langle x, e_i \rangle_< \in \Z$ because $x \in K_1^*$. We set 
\[
y:=\sum_{j \in S} a_j e^*_j + \sum_{j \notin S} \left\lceil a_j \right\rceil e_j^* \in K_2^*
\]
and note that 
\[
| x-y |_{\R^d} \leq \sum_{j \notin S} \left| a_j - \left\lceil a_j \right\rceil \right| | e_j^* |_{\R^d} \leq \sum_{j \notin S} | e_j^* |_{\R^d} <\eps,
\]
which proves Theorem \ref{cor::polarlattice}.
\end{proof}
Lastly,  we need to prove Proposition \ref{prop::boundedn0}.  For this, we use essentially a compactness argument, which relies on the lemma below.
\begin{lemma} \label{lem::opensubgp}
Let $S$ be a finite set of primes and let $V_p$ be a finite dimensional $\Q_p$-vector space with the topology induced by $\Q_p$ for all $p \in S$.  Set $V=\prod_{p \in S} V_p$ (equipped with the product topology) and let $H \subset V$ be a closed additive subgroup.  

Then firstly,  $H=\prod_{p \in S} H_p$,  where each $H_p \subset V_p$ is a closed subgroup.  Furthermore,  $H$ is open in $W:= \prod_{p \in S} W_p$, where $W_p \subset V_p$ is the smallest $\Q_p$-subspace containing $H_p$.
\end{lemma}
\begin{proof}
The first part is due to the Chinese remainder theorem. Given $y=(y_p)_{p \in S} \in H$, we want to show that also $(y_{p_1}, 0, \dots, 0) \in H$. To this end, we find for any $r \in \N$ an $m_r \in \Z$ such that $m_r \equiv \delta_{p=p_1} \mo p^r$ for all $p \in S$, which is possible by the CRT. Then $m_r y \to (y_{p_1}, 0, \dots, 0)$, which shows the claim because $H$ is closed.

For the second part,  it is enough to show that $H_p$ is open in $W_p$ for each $p \in S$. This is proved in Chapter 2 as Proposition 5  in \cite{weil}.
\end{proof}
\begin{proof}[Proof of Proposition \ref{prop::boundedn0}]
Let $v \in \bigcup_{n \geq 0} \lambda^n L^*$ be given.  If $v$ does not have property $(\star)$, we are done.  If it does, there is $y_1 \in U$ such that $\langle v, y_1 \rangle_< \in \frac{1}{M^s} \Z$ but $\langle v, y_1 \rangle_< \notin \Z$. We set $H_1:=\overline{\langle y_1 \rangle} \subset A_\beta$, where the closure is with respect to the topology of $A_\beta$. We recall that as an Abelian group,  $A_\beta=\prod_{p|M} V_p$ for $V_p$ finite dimensional $\Q_p$-vector spaces. Thus, by Lemma \ref{lem::opensubgp}, $H_1$ is of the form $H_1=\prod_{p|M} H_1^{(p)}$. We let
\[
W_1^{(p)}:= \left\langle H_1^{(p)} \right \rangle_{\Q_p},
\]
for $p|M$ be the subspaces generated by $H_1^{(p)}$ and let $W_1=\prod_{p|M} W_1^{(p)}$. As $\langle v, y_1 \rangle_< \notin \Z$,  $y_1$ is non-zero,  so at least one of them must have non-zero dimension.

By Lemma \ref{lem::opensubgp}, $M^sH_1$ is open in $W_1$.  As $L$ is bounded in $\Ab$ and $\beta$ is a contraction, there is $m_1^\prime=m_1^\prime(y_1)$ such that for all $m \geq m_1^\prime$,
\[
\beta^{m} L \cap W_1 \subset M^s H_1.
\]
We let 
\[
m_1:=\max_{y_1 \in U} m_1^\prime\left(y_1\right)
\]
which is well-defined as $U$ is finite and, crucially, independent of $v$.  Taking the minimum over $U$ there (and proving that this $U$ is finite in the first place)   is necessary to assure that $m_1$ is independent of $v$ and $n$.  We note that then for any $z \in \beta^{m_1} L \cap W_1$ one has that $\langle v, z \rangle_< \in \Z$ because $\langle v,  y \rangle_< \in M^{-s} \Z$ for any $y \in H_1$. 

For the second step, we consider $\beta^{m_1} v$.  Either $\beta^{m_1} v$ does not have property $(\star)$ or there is $y_2 \in U$ such that $\langle v, \beta^{m_1} y_2 \rangle_< \in \frac{1}{M^s} \Z, \langle v, \beta^{m_1} y_2 \rangle_< \notin \Z$.  As $\beta^{m_1} y_2 \in \beta^{m_1} L$,  we deduce from the last sentence in the preceding paragraph that $\beta^{m_1} y_2 \notin W_1$. We define 
\[
H_2:=\overline{\langle y_1, \beta^{m_1} y_2 \rangle}, 
\]
let $W_2$ be as above the smallest product of subspaces containing $H_2$ and note that $W_2 \supsetneq W_1$.
We once again apply Lemma \ref{lem::opensubgp} to find $m_2^\prime=m_2^\prime(y_1, y_2)$ such that for all $m \geq m_2^\prime$,
\[
\beta^{m} L \cap W_2 \subset M^s H_2.
\]
We set 
\[
m_2:=\max_{z_1, z_2 \in U} m_2^\prime(z_1, z_2)
\]
and note that then for all $z \in \beta^{m_2} L \cap W_2$, $\langle v, z \rangle_< \in \Z$, independently of $v$. 

We now iterate this procedure. As the dimension of $W_i$ increases at each step and $A_\beta$ is finite dimensional, this must terminate after finitely many steps, say after $r=r(d, |\{p: p|M\}|)$ many. We let 
\[
n_0:=m_r:=\max_{z_1, \dots , z_r \in U} m_r^\prime(z_1, \dots , z_r)
\]
be the maximum after the $r$-th step defined just as above. If $\beta^m v$ has property $(\star)$ for all $m<n_0$,  then $\beta^{n_0}$ cannot have property $(\star)$. This is because if it had property $(\star)$,  one could do another iteration of the procedure above - increasing the dimension,  which is impossible by definition of $r$.
\end{proof}

\printbibliography
\end{document}